\documentclass[article,onefignum,onetabnum]{siamart220329}

\usepackage{lipsum}
\usepackage{amsfonts}
\usepackage{graphicx}
\usepackage{epstopdf}
\usepackage{algpseudocode}
\ifpdf
  \DeclareGraphicsExtensions{.eps,.pdf,.png,.jpg}
\else
  \DeclareGraphicsExtensions{.eps}
\fi

\usepackage{mathtools}
\usepackage{bm}
\usepackage{pgfplots}
\usepackage{tensor}
\usetikzlibrary{pgfplots.groupplots}
\usepackage{multirow}
\usepackage[ntheorem]{empheq}
\usepackage{cases}
\usepackage[export]{adjustbox}
\usepackage[dvipsnames]{xcolor}
\usepackage{amsopn}

\newcommand{\cC}{\ensuremath{\mathcal{C}}}
\newcommand{\cD}{\ensuremath{\mathcal{D}}}
\newcommand{\cE}{\ensuremath{\mathcal{E}}}

\newcommand{\cH}{\ensuremath{\mathcal{H}}}

\newcommand{\cJ}{\ensuremath{\mathcal{J}}}

\newcommand{\cM}{\ensuremath{\mathcal{M}}}

\newcommand{\cP}{\ensuremath{\mathcal{P}}}

\newcommand{\bA}{\ensuremath{\mathbb{A}}}
\newcommand{\bB}{\ensuremath{\mathbb{B}}}
\newcommand{\bC}{\ensuremath{\mathbb{C}}}

\newcommand{\bG}{\ensuremath{\mathbb{G}}}

\newcommand{\bM}{\ensuremath{\mathbb{M}}}
\newcommand{\bN}{\ensuremath{\mathbb{N}}}

\newcommand{\bR}{\ensuremath{\mathbb{R}}}
\newcommand{\bS}{\ensuremath{\mathbb{S}}}

\newcommand{\bV}{\ensuremath{\mathbb{V}}}

\newcommand{\Ccal}{\ensuremath{\mathcal{C}}}

\newcommand{\Ecal}{\ensuremath{\mathcal{E}}}

\newcommand{\Hcal}{\ensuremath{\mathcal{H}}}

\newcommand{\Jcal}{\ensuremath{\mathcal{J}}}

\newcommand{\Vcal}{\ensuremath{\mathcal{V}}}

\newcommand{\Vn}{\ensuremath{V_n}}
\newcommand{\Wm}{\ensuremath{W_m}}
\newcommand{\Vtwon}{\ensuremath{V_{2n}}}
\newcommand{\Wtwom}{\ensuremath{W_{2m}}}

\newcommand{\ms}{\ensuremath{\mathrm{ms}}}

\newcommand{\prm}{\ensuremath{\theta}}

\newcommand{\rd}{\ensuremath{\mathrm d}}

\newcommand{\dt}{\ensuremath{\mathrm dt}}

\DeclareMathOperator*{\vspan}{span}
\newcommand{\rank}{\ensuremath{\text{rank}}}
\newcommand{\eps}{\ensuremath{\varepsilon}}

\newcommand{\Nr}{\ensuremath{{2n}}} 
\newcommand{\Nrh}{\ensuremath{n}} 

\newcommand{\R}[2]{\ensuremath{\mathbb{R}^{{#1}\times{#2}}}}

\newcommand{\J}[1]{\ensuremath{J_{#1}}} 
\newcommand{\Idm}{\ensuremath{I}} 
\newcommand{\Zrm}{\ensuremath{0}} 

\newcommand{\norm}[1]{\|{#1}\|}

\newcommand\dist{\text{dist}}
\newcommand{\ord}{\ensuremath{\mathcal{O}}}

\newcommand{\cond}{\; :\;}

\DeclareMathOperator{\argmin}{argmin}
\newcommand{\Vr}{V_{\Nr}}
\newcommand{\Mr}{M_{\Nr}}
\newcommand{\ur}{u_{\Nr}}
\newcommand{\dur}{\dot{u}_{\Nr}}

\newcommand{\uq}{u^q}
\newcommand{\up}{u^p}
\newcommand{\inprod}[3]{\left<{#1},{#2}\right>_{#3}}
\newcommand{\inprodV}[2]{\ensuremath{\left<{#1},{#2}\right>}}
\newcommand{\Vqp}{\ensuremath{\widehat{V}}} 

\newcommand{\xs}{{\overline{x}}} 
\newcommand{\Sdom}{\ensuremath{D}} 

\newcommand{\Thtrain}{\Theta_h}
\newcommand{\Thtest}{\Theta_s}

\newcommand{\Ese}{\Ecal}
\newcommand{\Eproj}{\underline{\Ese}}
\newcommand{\Ebound}{\overline{\Ese}}

\newcommand{\Esets}{\Ese_{\Thtest}}
\newcommand{\Esetsstat}{\Esets^{st}}
\newcommand{\Esetsdyn}{\Esets^{dyn}}
\newcommand{\Eprojts}{\underline{\Ese}_{\Thtest}}
\newcommand{\Eboundts}{\overline{\Ese}_{\Thtest}}

\newcommand{\eHts}{\Ecal_{\Hcal,\Thtest}}
\newcommand{\dH}{\Delta\Hcal}
\newcommand{\dHts}{\dH_{\Thtest}} 

\newsiamremark{remark}{Remark}
\newsiamremark{hypothesis}{Hypothesis}
\crefname{hypothesis}{Hypothesis}{Hypotheses}
\newsiamthm{claim}{Claim}

\headers{Dynamical PBDW for Filtering Problems}{O. Mula, C. Pagliantini, and F. Vismara}

\title{
Dynamical approximation and sensor placement \\for filtering problems}

\author{Olga Mula\thanks{Centre for Analysis, Scientific computing and Applications, Eindhoven University of Technology, The Netherlands. (\email{o.mula@tue.nl}, \email{f.vismara@tue.nl}).}\and Cecilia Pagliantini\thanks{Dipartimento di Matematica, Universit\`a di Pisa, Italy. (\email{cecilia.pagliantini@unipi.it}).}\and Federico Vismara\footnotemark[1]}

\pgfplotsset{compat=1.18}
\begin{document}

\maketitle

\begin{abstract}
We consider the inverse problem of reconstructing an unknown function $u$ from a finite set of measurements, under the assumption that $u$ is the trajectory of a transport-dominated problem with unknown input parameters.
We propose an algorithm based on the Parameterized Background Data-Weak method (PBDW) where dynamical sensor placement is combined with approximation spaces that evolve in time. We prove that the method ensures an accurate reconstruction at all times and allows to incorporate relevant physical properties in the reconstructed solutions by suitably evolving the dynamical approximation space.
As an application of this strategy we consider Hamiltonian systems modeling wave-type phenomena, where preservation of the geometric structure of the flow plays a crucial role in the accuracy and stability of the reconstructed trajectory.
\end{abstract}

\section{Introduction}
\subsection{The filtering problem}
\label{sec:intro-filtering}
Data assimilation of problems presenting physical properties such as conservation laws, discontinuities and strong transport effects is still nowadays a challenging task despite the long history of the field. In this work, we address this topic by developing a strategy for online filtering which is especially thought for this type of dynamics.

In the following, we define filtering as the following task. Let $\Sdom\subseteq \bR^d$ be a spatial domain and let $V$ be a Hilbert space with norm $\norm{\cdot}$ induced by the inner product $\left< \cdot, \cdot \right>$.
We consider a time-dependent process in which, at every time $t\geq 0$, the goal is to recover the unknown state $u(t)\in V$ from data given by $m$ partial observations
\begin{equation*}
    z_i(t) = \ell_i(u(t)), \qquad i = 1,\dots, m,
\end{equation*}
where $\ell_i: V \rightarrow \bR$ are known, independent, continuous linear functionals from $V'$, the dual of $V$. The functionals $\ell_i$ often represent the response of a physical measurement device but they can have a different interpretation depending on the application. The knowledge of $z(t)=(z_i(t))_{i=1}^m$ is equivalent to that of the orthogonal projection $\omega(t)=P_{\Wm}u(t)$, where
\begin{equation*}
    \Wm = \vspan\{ \omega_1, \dots, \omega_m \},
\end{equation*}
and $\omega_i\in V$ are the Riesz representers of the linear functionals $\ell_i$, that is,
\begin{equation*}
    \ell_i(v) = \left< \omega_i, v\right> \qquad \forall\, v\in V.
\end{equation*}
Recovering $u(t)\in V$ from $\omega(t)\in \Wm$ is ill-posed as soon as the dimension of $V$ exceeds $m$. Indeed, in that case, for any observation $\omega(t)\in \Wm$, there are infinitely many $v(t)\in V$ such that $P_{\Wm}v(t)=\omega(t)$. Thus the only way to recover the evolution of $u$ up to a guaranteed accuracy is to combine the measurements with some a priori information. In our case, $u$ is assumed to be the solution to an evolution problem of the form
\begin{equation}\label{eq:PDE}
    \begin{cases}
        \partial_t u - \cP_\theta(u) = 0 & \quad x\in \Sdom,\, t\in\bR_+, \\
        u(t=0, x, \theta) = u_0(\theta) & \quad x\in\Sdom,
    \end{cases}
\end{equation}
where $\cP_\theta$ is a partial differential operator (with appropriate boundary conditions). The operator depends on a parameter $\theta$ from a compact space $\Theta$, which can be finite or infinite dimensional. We assume that, for each $\theta \in \Theta$, the problem is well-posed in the sense that there exists a unique solution $u(\theta)\in \cC(\bR_+; V)$. This notation means that, for each $\theta\in \Theta$, $u(\theta)$ is a function of $t$ and $x$; the evolution is continuous in time, and for each $t\in \bR_+$, $u(\theta)(t, \cdot)\in V$. To simplify notation, we will sometimes write $u(\theta)(t, x)$ as $u(t, x, \theta)$.

In the context of the filtering problem, the parameter $\theta^\dagger$ corresponding to the measurements is unknown, so the dynamics of $u(\theta^\dagger)$ cannot be recovered by solving problem \eqref{eq:PDE} with direct methods. Instead, we are given the partial observations $z(t)$, and our prior assumption on $u$ is that it belongs to the set of trajectories
\begin{equation*}
\cM \coloneqq \{ u(\theta) \in \cC(\bR_+; V) \cond \theta \in \Theta  \}.
\end{equation*}
For every $t\in \bR_+$, the goal is to recover $u(\theta^\dagger)(t,\cdot)$ using the observations $z(\tau)$ for $\tau \in [0, t]$, and the fact that $u(\theta^\dagger)\in \cM$.


\subsection{Dynamical approximation of the reconstruction space}
One natural strategy to address the complex structure of $\cM$ is to search for a recovery of $u(\theta^\dagger)(t, \cdot)$ by an element from a low-dimensional reconstruction space $\Vn\subset V$. The space $\Vn$ could be either an $n$-dimensional linear subspace or, more generally, a nonlinear approximation space parameterized by $n$ degrees of freedom, with $n\leq m$. In order to obtain quantitative results in the recovery, $\Vn$ needs to approximate the elements of $\cM$ as well as possible in the sense that
\begin{equation*}
    \dist(\cM, \Vn) \coloneqq \max_{\theta\in \Theta}\, \max_{t\in \bR_+}\,\min_{v\in \Vn} \Vert u(\theta)(t, \cdot) - v \Vert
\end{equation*}
should decay fast with $n$, and we need that the quantity is small already for moderate values of $n$. When $\Vn$ is restricted to be a linear subspace, the ideal benchmark is given by the Kolmogorov $n$-width
\begin{equation}\label{eq:kolmo-static}
    d_n(\cM) \coloneqq \inf_{\dim(Y_n)\leq n} \dist(\cM, Y_n),
\end{equation}
where the infimum is taken over all linear subspaces $Y_n$ of dimension at most $n$. If $\cM$ encodes a family of trajectories that develop complex features in time, $d_n(\cM)$ will decay slowly with $n$. This is expected to happen for several types of dynamics, in particular for systems where the support of the solution is transported in time. We illustrate this claim through a simple example later on. In these cases there will not be a $\Vn$ that accurately approximates $u(\theta)(t, \cdot)$ uniformly at all times and for all parameters. This naturally leads us to slice the set $\cM$ in time as
\begin{equation*}
    \cM = \cup_{t\geq 0} \cM(t),
\end{equation*}
with
\begin{equation}\label{eq:solution-set-t}
    \cM(t) \coloneqq \{ u(\theta)(t,\cdot)\in V \cond \theta \in \Theta \},
\end{equation}
and consider dynamical reconstruction spaces $\Vn(t)$. In this setting, $\Vn(t)$ should be such that
\begin{equation*}
    \dist(\cM(t), \Vn(t)) \coloneqq \max_{\theta\in \Theta} \min_{v\in \Vn(t)} \Vert u(\theta)(t, \cdot) - v \Vert
\end{equation*}
decays fast as $n\to \infty$ and for all times $t\geq0$. For a given time slice $\cM(t)$, the best approximation in terms of linear subspaces is now given by the Kolmogorov $n$-width of $\cM(t)$ instead of $\cM$,
\begin{equation}\label{eq:kolmo-dynamic}
    d_n(\cM(t)) = \inf_{\dim(Y_n)\leq n} \dist(\cM(t), Y_n).
\end{equation}

Given the above conceptual ingredients, the rationale behind working with a dynamical $\Vn(t)$ is justified by the fact that we not only have that
\begin{equation*}
    d_n(\cM(t)) \leq d_n(\cM)\qquad \forall\, t\in\bR_+,
\end{equation*}
but we also hope to be in a regime where $d_n(\cM(t))$ is significantly smaller than $d_n(\cM)$. This gives room to build spaces $\Vn(t)$ with much higher approximation accuracy than the ones we would obtain with a constant $\Vn$. In other words, if $d_n(\cM(t)) \ll d_n(\cM)$, we expect to find practical algorithms to build $\Vn(t)$ such that $\dist(\cM(t), \Vn(t)) \ll \dist(\cM, \Vn)$ for all $t\geq 0$.

Apart from very simplified cases, both static and dynamic reconstructions spaces, $\Vn$ or $\Vn(t)$, achieving the above infima \eqref{eq:kolmo-static}-\eqref{eq:kolmo-dynamic} are usually out of reach. For static spaces, practical model reduction techniques such as polynomial approximation in the parameterized domain \cite{CD2015acta,CDS2011} or reduced bases \cite{EPR2010, HRS2015, RHP2007, ZKA2019} construct spaces $\Vn$ that are ``suboptimal yet good''. In particular, the reduced basis method, which constructs $\Vn$ as the span of a specific selection of particular solution instances $u^1,\dots,u^n\in \cM$, has been proven to have approximation error $\dist(\cM,\Vn)$ that decays with the same polynomial or exponential rates as $d_n(\cM)$, and, in that sense, is close to optimal \cite{DPW2013}. In the case of time-dependent problems, the use of low-dimensional time-dependent approximation spaces $\{\Vn(t)\}_{t\geq 0}$ is currently attracting a lot of attention for forward reduced modeling (see, e.g., \cite[Section 4.3]{HPR22} and references therein) but there are very few works leveraging this concept for data assimilation tasks. To the best of our knowledge, only \cite{Lombardi2022, SHNT2023, DY22} and \cite[Chapter 7]{Vidlickova2022} have considered dynamical approximation techniques for linear Kalman filtering.
In our work we propose an adaptive strategy for the construction and update of $\Vn(t)$ based on dynamical low-rank approximation \cite{KL07} and dynamically orthogonal (DO) schemes \cite{SL09}.

\subsection{Contributions and outline of the present work}
There are three main contributions in this work. The first and most important one is the development of an algorithm for online filtering that leverages the advantages of dynamical approximation spaces $\Vn(t)$, and for which we prove rigorous recovery bounds. It is especially suited for problems transporting sharp gradients, and localized solutions. Such problems arise in many different applications such as fluid dynamics or wave propagation to name a few.

Our strategy takes the Parameterized Background Data-Weak algorithm (PBDW) as a starting point. PBDW was originally introduced in \cite{MPPY2015} and it has been analyzed and extended in a series of papers such as \cite{BCDDPW2017, CDDFMN2020, CDMN2022, CDMS2022}. It has already found numerous applications ranging from biomedical problems \cite{TPYP2018, GLM2021, GLM2022} to pollution dispersion \cite{HCBM2019} and nuclear engineering \cite{ABCGMM2018, TLMMT2023}. We refer to \cite{Mula2022} for an overview. In its vanilla version, PBDW involves static, linear spaces $\Vn$. We develop an extension to dynamical spaces $\Vn(t)$, which we call Dynamical PBDW (Dyn-PBDW). We justify the potential gain of working with $\Vn(t)$ by building, in \Cref{sec:why-dynamical}, a simple, manufactured set $\cM$ that mimics transport dominated problems. For this set, $d_n(\cM(t))$ decays exponentially fast in the dimension $n$, while $d_n(\cM)$ decays as bad as $n^{-1/2}$.

As an additional challenge, transport-dominated problems are characterized by solutions that are often very localized in space, and propagate fronts at finite speed. If the sensors have fixed locations in space and capture only local features, then the support of the solution will be at times out of the range of the sensors. We show that such a situation leads to instabilities, and even ill-posedness of our recovery algorithm. It is thus necessary to move the sensors in such a way that they follow the dynamics in order to stably recover the solution. The second contribution of the paper is to provide an analytical framework to study this issue, and to provide an algorithm to dynamically move the sensors. Our strategy is based on maximizing the stability of the recovery algorithm in the sense that we define in Section \ref{sec:dyn_sens_plac}. Note that, with respect to the above notation, moving the sensors is equivalent to considering time-dependent observation spaces $\Wm(t)$. In essence, the strategy that we present leads to an evolution equation on the location of the sensors, which in turn defines an evolution equation for $\Wm(t)$.

Our third and last contribution is the application of the method to Hamiltonian systems. This family of problems is an excellent source of examples where filtering is challenged by both strong advection effects and the presence of a rich geometric structure (symplecticity of the phase space and conservation of invariants). To the best of our knowledge, our work provides the first filtering algorithm especially adapted to Hamiltonian problems, in the sense that the reconstructed solution belongs to a symplectic space and its Hamiltonian is preserved up to the reconstruction error.

The paper is organized as follows. In Section \ref{sec:earlier-works} we outline the connections of our contribution with previous works. In Section \ref{sec:DPBDW} we present Dyn-PBDW, and theoretical arguments in favor of working with dynamical spaces $V_n(t)$ and $W_m(t)$. Section \ref{sec:dyn_sens_plac} is devoted to a strategy for dynamical sensor placement based on maximizing the stability of our reconstruction method. Section \ref{sec:Ham_dynamics} explains how to apply the method to Hamiltonian systems with symplectic spaces $V_n(t)$. In Section \ref{sec:num_exp} we illustrate the performances of the Dyn-PBDW algorithm in 1D and 2D Hamiltonian systems involving the Schr\"odinger equation and the shallow water equations.

\subsection{Connections with earlier works}\label{sec:earlier-works}
Our analysis and algorithms share connections with earlier works which we next comment on.
\begin{itemize}
    \item $\Vn(t)$: Using dynamical approximation spaces $\Vn(t)$ for filtering is a rather unexplored idea. As already brought up, to the best of our knowledge only \cite{Lombardi2022, SHNT2023, DY22} and \cite[Chapter 7]{Vidlickova2022} have considered such techniques for linear Kalman filtering. Our approach significantly differs from these works in the sense that we do not follow a Bayesian approach. We rely on a more deterministic point of view where we do not assume any prior distribution over the parameters space. This results in more deterministic notions of accuracy quantification.
    \item $\Wm(t)$: Dynamically moving sensors is already done in certain applications (see, e.g., \cite{Haik23}) but this is the first time that a dynamical placement strategy is built on the basis of rigorous theoretical arguments of stability of the reconstruction. The key for this is to view the task of moving sensors as the one of making an observation space evolve in time. This point of view is novel, and we expect that it will inspire future works.
    \item Filtering of Hamiltonian dynamics has never been the subject of a dedicated study to the best of our knowledge. Interestingly, this does not mean that Hamiltonian dynamics are absent in the literature of data assimilation: Langevin and Hamiltonian-based Monte Carlo algorithms for Bayesian filtering are actually a topic of very active research (see, e.g., \cite{GC2011}). The way Hamiltonian dynamics arises in our work is therefore totally different.
    \item PBDW: The interplay between reduced order modeling and inverse problems has been the subject of very active research in the past years, and PBDW has played a central role in this trend. From the point of view of this body of literature, our main contribution  lies first of all in the type of inverse problem that we address, namely filtering. Most existing works concentrate either on problems that do not depend on time or that are focused on the so-called smoothing task for time-dependent problems. Smoothing consists in recovering the state $u$ in a whole time interval $[0, T]$ from all the observations in that time frame. It is different from filtering where one recovers the state $u(t)$ online as time $t$ advances. Since there was no strategy to update the spaces $\Vn(t)$ and $\Wm(t)$ in earlier works, addressing filtering and transport dominated problem was not really possible. It is however worth mentioning that the topic of optimal sensor placement in connection to PBDW and close relatives like the Generalized Empirical Interpolation Method has been extensively studied in, e.g., \cite{MM2013, MMPY2015,BCMN2018}. In these works, the sensors locations are chosen one after the other in a greedy manner. Beyond the dynamical component of our algorithm, our approach updates all locations in a joint optimization process.
\end{itemize}

\section{A dynamical PBDW for filtering problems}\label{sec:DPBDW}
\subsection{The PBDW method}\label{sec:DPBDW-method}
In its vanilla version, a reconstruction with the PBDW method involves a static, linear observation space $\Wm$ and a static, linear reduced space $\Vn$. The reconstruction is defined as the mapping $A: \Wm\to \Vn\oplus(\Wm\cap\Vn^\perp)$ such that
\begin{equation*}
    A(\omega) = \argmin_{v\in \omega +\Wm^\perp} \Vert v - P_{\Vn} v \Vert \qquad \forall\, \omega \in \Wm.
\end{equation*}
The extension to dynamical $\Vn(t)$ and $\Wm(t)$ is therefore straightforward. For all $t\in \bR_+$, we define the mapping $A(t): \Wm(t) \to \Vn(t)\oplus(\Wm(t)\cap\Vn^\perp(t))$ such that
\begin{equation*}
    A(t)(\omega(t)) = \argmin_{v\in \omega(t) +\Wm^\perp(t)} \Vert v - P_{\Vn(t)} v \Vert \qquad \forall\, \omega(t) \in \Wm(t).
\end{equation*}
Denoting
\begin{equation}\label{eq:beta}
    \beta(X, Y) \coloneqq \inf_{x\in X} \sup_{y\in Y} \dfrac{\left< x, y \right>}{\Vert x \Vert \Vert y\Vert} = \inf_{x\in X} \dfrac{\Vert P_Y x \Vert}{\Vert x \Vert} \; \in [0, 1]
\end{equation}
for any pair of subspaces $(X, Y)$ of $V$, the above optimization problem has a unique minimizer which can be written as
\begin{equation}\label{eq:At}
    A(t)(\omega(t)) = u^*_t(\omega(t)) = v^*_t(\omega(t)) + \omega(t) - P_{\Wm(t)} v^*_t(\omega(t))
\end{equation}
provided that
\begin{equation}\label{eq:cond-pbdw}
    1\leq n\leq m \quad\text{ and }\quad \beta(\Vn(t), \Wm(t))>0.
\end{equation}
We adhere to these assumptions in the following and we refer to, e.g., \cite[Lemma A.1]{Mula2022} for the proof of this result. The element $v^*_t(\omega(t))$ is defined as
\begin{equation}\label{eq:vtstar}
    v^*_t(\omega(t)) \coloneqq \left( P_{\Vn(t)|\Wm(t)} P_{\Wm(t)|\Vn(t)} \right)^{-1} P_{\Vn(t)|\Wm(t)}(\omega(t)) \in \Vn(t),
\end{equation}
where, for any pair of closed subspaces $(X, Y)$ of $V$, $P_{X|Y}: Y \to X$ is the orthogonal projection onto $X$ restricted to $Y$. The operator $P_{\Vn(t)|\Wm(t)} P_{\Wm(t)|\Vn(t)}$ is invertible under the conditions \eqref{eq:cond-pbdw}.

It follows from formula \eqref{eq:At} that $A(t)$ is a bounded linear map from $\Wm(t)$ to $\Vn(t)\oplus (\Wm(t) \cap \Vn^\perp(t))$. It was proven in \cite[Theorem 2.9]{BCDDPW2017} that the reconstruction error satisfies the bound
\begin{equation}
    \label{eq:err-pbdw-ustar}
    \Vert u(t) - u^*_t(\omega(t))\Vert \leq\dfrac{1}{\beta(\Vn(t), \Wm(t))} \eps(u(t), V_n(t), \Wm(t)) \qquad \forall\, u(t) \in V,
\end{equation}
where
\begin{equation*}
    \eps(u(t), V_n(t), \Wm(t)) \coloneqq \inf_{v(t)\in \Vn(t)\oplus (\Wm(t) \cap \Vn^\perp(t))} \Vert u(t)-v(t)\Vert.
\end{equation*}
We quantify the performance of the algorithm over all possible states $u(t)\in\cM(t)$ as
\begin{equation}
\label{eq:performance-wc}
E(A(t), \cM(t)) \coloneqq \max_{u(t)\in \cM(t)} \Vert u(t) - A(t)(\omega(t))\Vert.
\end{equation}
Thanks to \eqref{eq:err-pbdw-ustar}, we have that
\begin{align}
\begin{split}\label{eq:bound-u-star}
    E(A(t), \cM(t)) &\leq \dfrac{\dist(\cM(t), \Vn(t) \oplus (\Wm(t) \cap \Vn^\perp(t)))}{\beta(\Vn(t), \Wm(t))} \leq \dfrac{ \dist(\cM(t), \Vn(t))}{\beta(\Vn(t), \Wm(t))}.
\end{split}
\end{align}
Note that the performance of the static PBDW, $E(A, \cM)$, satisfies the same relation but with no time-dependency in the above expression.

We now pause the main discussion to make three remarks. First, it is important to see that $\beta(\Vn(t), \Wm(t))$ plays the role of a stability constant. It can be interpreted as the cosine of the angle between the spaces $\Vn(t)$ and $\Wm(t)$. To have a stable reconstruction, we need to guarantee that $\beta(\Vn(t), \Wm(t))\geq \underline{\beta}$ for a $\underline{\beta}>0$ and for all $t\geq 0$. Ideally, we would like to have $\underline{\beta}=1$ but this would require that $\Vn(t)\subseteq \Wm(t)$. Satisfying this inclusion is not possible in practice because the physics of the sensors makes the observation functions $\omega_i(t)$ very different in nature compared to the elements of the space $V_n(t)$ (which encode the physics of the problem).

The second remark is about the error bound \eqref{eq:bound-u-star}. Note that the right-hand side involves the distance of $u(t)$ to $\Vn(t) \oplus (\Wm(t) \cap \Vn^\perp(t))$. Since in general $\Wm(t) \cap \Vn^\perp(t)$ provides very little approximation power, it usually suffices to work with $v^*_t(\omega(t))$, whose reconstruction performance satisfies
\begin{equation}\label{eq:Evtstar}
    E(v_t^*,\cM(t)):=\max_{u(t)\in \cM(t)} \Vert u(t) - v^*_t(\omega(t))\Vert \leq \dfrac{\dist(\cM(t), \Vn(t))}{\beta(\Vn(t), \Wm(t))} .
\end{equation}
This is interesting when $\Vn(t)$ is built to preserve geometrical properties of the original problem. In our application, we work with $v^*_t\in\Vn(t)$ instead of $u^*_t$ for this reason.

Our last remark is about the fact that our performance criterion \eqref{eq:performance-wc} is defined in the worst case sense: the element of $\cM(t)$ which is reconstructed at worst determines $E(A(t), \cM(t))$. This is a major distinction with respect to the Bayesian approach to filtering which assumes a probability distribution $\nu_V(t)$ on $V$ which is supported over $\cM(t)$, and which usually evolves in time. The distribution is typically the result of push-forwarding a probability distribution $\nu_\Theta(t)$ on the parameter set $\Theta$ through the parameter-to-solution map: $\nu_V(t) = u \# \nu_\Theta(t)$. In the Bayesian approach, one typically starts from  a given distribution $\nu_\Theta(0)$ at the initial time, and $\nu_\Theta(t)$ evolves following different possible update strategies. The performance is then measured in an average sense, usually in terms of the mean-square error
$$
E_{\ms}(A(t), \cM(t), \nu(t)) \coloneqq
\left( \int_V \Vert u(t) - v^*_t(\omega(t))\Vert^2 \rd\nu_V(u) \right)^{1/2}.
$$
We may observe that, for any prior distribution $\nu_\Theta(t)$,
$$
E_{\ms}(A(t), \cM(t), \nu(t))
\leq
E(A(t), \cM(t), \nu(t)),
$$
thus the worst-case performance criterion could be understood as a distributionally robust way of measuring performance.

\subsection{Why dynamical \texorpdfstring{$\Vn(t)$}{Vnt} and \texorpdfstring{$\Wm(t)$}{Wmt}?}\label{sec:why-dynamical}
The following simple manipulation provides an argument to justify why working with dynamical $\Vn(t)$ and $\Wm(t)$ can give a very significant gain in reconstruction accuracy. Since $E(A, \cM) \geq \dist(\cM, \Vn)$, we have that
\begin{equation}\label{eq:ratio-static-dyn-1}
    \frac{E(A(t), \cM(t))}{E(A, \cM)} \leq \dfrac{1}{\beta(\Vn(t), \Wm(t))} \frac{\dist(\cM(t), \Vn(t))}{\dist(\cM, \Vn)}.
\end{equation}
Now, assuming that we can build a near optimal $\Vn(t)$ in the sense that there exists a $C\geq 1$ of moderate value such that $\dist(\cM(t), \Vn(t)) \leq C d_n(\cM(t))$, we can further bound \eqref{eq:ratio-static-dyn-1} as
\begin{equation*}
    \frac{E(A(t), \cM(t))}{E(A, \cM)} \leq \dfrac{C}{\beta(\Vn(t), \Wm(t))} \frac{d_n(\cM(t))}{d_n(\cM)}.
\end{equation*}
If we place the sensors in such a way that $\Wm(t)$ yields $\beta(\Vn(t), \Wm(t))\geq \underline{\beta}$ for some $\underline{\beta}>0$, then
\begin{equation*}
    \frac{E(A(t),\cM(t))}{E(A, \cM)} \leq \dfrac{C}{\underline{\beta}} \frac{d_n(\cM(t))}{d_n(\cM)}, \quad \forall\, t\geq 0.
\end{equation*}
As we anticipated in the introduction, having $\frac{d_n(\cM(t))}{d_n(\cM)} \ll 1$ leads to an important gain in reconstruction accuracy if the sensors are placed in a way that guarantees enough stability.

The following simple example justifies the fact that, in problems with strong advection effects, we can indeed expect that $\frac{d_n(\cM(t))}{d_n(\cM)} \ll 1$ for all $t\geq0$. Consider the pure transport equation in $V=L^2(\Sdom)$ with $\Sdom=\bR$. The problem involves the parameter $\theta = (\theta_1, \theta_2)$ in a compact subset $\Theta$ of $\bR^2$, and reads
\begin{equation*}
\begin{cases}
    \partial_t u(t, x, \theta) + \theta_2 \partial_x u(t, x, \theta)=0,  \\ u(0, x, \theta)= \dfrac{1}{\sqrt{2\pi}\,\theta_1}  \exp{\left(-\dfrac{ x^2 }{2\theta_1^2}\right)}.
\end{cases}
\end{equation*}
The solution to this PDE is known analytically, and reads
\begin{equation*}
    u(t, x, \theta) = \dfrac{1}{\sqrt{2\pi}\,\theta_1}  \exp{\left(-\dfrac{ (x-t\theta_2)^2 }{2\theta_1^2}\right)}.
\end{equation*}
The set of parametric solutions is then
$
    \cM = \{ u(\theta) \cond \theta \in \Theta \} \subset \cC(\bR_+; L^2(\bR))
$
and its time-sliced version is
\begin{equation*}
    \cM(t) = \left\{ x\mapsto \dfrac{1}{\sqrt{2\pi}\,\theta_1}  \exp{\left(-\dfrac{ (x-t\theta_2)^2 }{2\theta_1^2}\right)} \cond \theta \in \Theta \right\}  \subset L^2(\bR).
\end{equation*}
Since, for a given $t\geq 0$, $\theta\in \Theta \mapsto u( \theta)(t, \cdot) \in L^2(\bR)$ is analytic in $\theta$, $d_n(\cM(t))$ decays exponentially in $n$. However, since the time interval $\{t\geq 0\}$ is not bounded, $d_n(\cM) = \ord(1)$ for all $n\in \bN$. One may wonder to what extent the behavior of $d_n(\cM)$ is due to the fact that the time interval is unbounded. In fact, if we consider a bounded interval $[0, T]$, it is possible to build examples for which we have $d_n(\cM)= \ord(1)$ for all $n\leq n_0$ and where $n_0$ depends on $T$ and $\Theta$.

The above example also gives a justification as to why we need to dynamically move the sensors. Suppose that $\Vn(t)= \vspan\{ u(\theta_i)(t, \cdot) \}_{i=1}^n$  for well chosen $\theta_i \in \Theta$ for $i=1,\dots, n$. Suppose further that the observation space is spanned by $m$ functions from the dictionary
\begin{equation*}
    \cD = \left\lbrace \omega_q \in V \;\cond \, \omega_q(x) = \dfrac{1}{\sqrt{2\pi}\,\sigma}  \exp{\left(-\dfrac{ (x-q)^2 }{2\sigma^2}\right)},\; \forall\, x\in D\right\rbrace \subset V,
\end{equation*}
where $\sigma>0$ is a fixed constant. If $\Wm= \vspan\{\omega_{q_i}\}_{i=1}^m$ for fixed $q_i\in \bR$, then $\left<\omega_{q_i}, u(\theta_j)(t, \cdot) \right>\to 0$  as $t\to \infty$ for all $1\leq i\leq m$ and $1\leq j\leq n$. As a result, $\beta(\Vn(t), \Wm)\to 0$ as $t\to \infty$.

\subsection{Noise and model errors}
For simplicity of exposition, we have omitted adding observation noise and modeling errors in the above analysis, but these sources of error can easily be included as follows. Let us assume that we get noisy observations $\widetilde{\omega} (t) = \omega(t) + \eta(t)$ with $||\eta(t)|| \leq \eps_{noise}(t)$. Suppose also that the true state $u(t)$ does not lie in $\cM(t)$ but satisfies $\dist(u(t), \cM(t)) \leq \eps_{model}(t)$. We can prove that the error bound \eqref{eq:err-pbdw-ustar} should be modified into
\begin{align*}
\Vert u(t) - u^*_t(\widetilde{\omega}(t)) \Vert
&\leq \frac{\eps(u(t), V_n(t), W_m(t)) + \eps_{noise}(t) + \eps_{model}(t)}{\beta(\Vn(t), \Wm(t))},
\end{align*}
and \eqref{eq:bound-u-star} becomes
\begin{align*}
E(A(t), \cM(t)) &\leq
\dfrac{\dist(\cM(t), \Vn(t) \oplus (\Wm(t) \cap \Vn^\perp(t))) + \eps_{noise}(t) + \eps_{model}(t)}{\beta(\Vn(t), \Wm(t))}.
\end{align*}
As one would expect, these estimates show that the reconstruction improves  when decreasing the model error and the noise.
Another source of inaccuracy comes from the modeling of the observation functions $\omega_i(t)$, for any $1\leq i\leq m$. Suppose that we work with imperfect functions $\widetilde{\omega}_i(t)$ that deviate from the exact ones $\omega_i(t)$ by $\Vert \omega_i(t) - \widetilde{\omega}_i(t) \Vert \leq \rho$ uniformly in $t\in \bR_+$ for some $\rho>0$. Then, noiseless observations can be written as
$
z_i(t) = \ell_i(u(t)) = \left< \omega_i(t), u(t) \right> = \left< \widetilde{\omega}_i(t), u(t)\right> + \left< \omega_i(t) - \widetilde{\omega}_i(t), u(t)\right>
$, for any $1\leq i\leq m$.
The right hand side tells us that, by working with the inexact $\widetilde{\omega}_i(t)$, we are introducing a noise term, which is $ \left< \omega_i(t) - \widetilde{\omega}_i(t), u(t)\right>\leq \rho \Vert u(t) \Vert$. This means that working with an inexact representation of the sensor response can be understood as introducing additional noise of level $\eps_{\text{noise}}(t)=\rho \Vert u(t) \Vert$ to the observations.

\section{Dynamical sensor placement}\label{sec:dyn_sens_plac}
The examples discussed in the previous section motivate the development of dynamical strategies to evolve both the approximation space $\Vn$ and the observation space $\Wm$ over time. The dynamical construction of the approximation space $\Vn(t)$ will be discussed in \Cref{sec:Ham_dynamics}, and, for the moment, we assume it to be prescribed at each time $t\in\bR_+$. We thus set out to define the evolution of the observation space $\Wm(t)$. The main idea is to use the stability constant $\beta(t)=\beta(\Vn(t),\Wm(t))$ from \eqref{eq:beta} to assess the quality of the reconstruction and drive the adaptation of $\Wm$ to maximize $\beta$.

Assuming that the $i$th measurement $\ell_i$ is a function depending on the position of the $i$th sensor, our goal is to maximize the stability constant $\beta$ -- thought of as a function of the observation space -- by dynamically adapting the sensor locations.
More in detail, we denote the position of the $i$th sensor, for $i=1,\dots,m$, at time $t$ by
\begin{equation*}
    \xs_i(t)=(\xs_i^{(1)}(t),\dots,\xs_i^{(d)}(t))\in\bR^d,
\end{equation*}
and we arrange the coordinates of the $m$ sensors in the vector
\begin{equation*}
    \xs(t)=(\xs^{(1)}(t),\dots,\xs^{(d)}(t))\in(\bR^m)^d\simeq\bR^{md}
\end{equation*}
where $\xs^{(\ell)}(t):=(\xs_1^{(\ell)}(t),\dots,\xs_m^{(\ell)}(t))\in\bR^m$, for $\ell=1,\dots,d$, collects the $\ell$th spatial component of all sensors. Then
\begin{equation*}
    \Wm(t)=\text{span}\{\omega_1(\xs_1(t)),\dots,\omega_m(\xs_m(t))\},
\end{equation*}
so that we can see the observation space as a function of $\xs(t)$, i.e., $\Wm=\Wm(\xs(t))$. Our task is to determine the sensor placement $\xs(t)$ at every time instant $t$. The strategy we propose is based on defining the optimal evolution of the sensors positions in the temporal interval $(0,T]$ as
\begin{equation}\label{eq:lagrangian}
    \xs = \mathop{\arg \max}\limits_{\substack{y : [0, T] \to \bR^{md} \\ y(0) = a \\ \dot{y}(0)=b}}\int_0^T \left( \beta^2\big(\Vn(\tau), \Wm(y(\tau))\big) - \frac{\lambda}{2} \Vert \dot{y}(\tau)\Vert^2 \right) \,d\tau,
\end{equation}
where $a$ and $b$ are suitably prescribed initial position and velocity of the sensors, respectively, while $\lambda\in\bR_+$ is the coefficient of a penalization term proportional to the sensors velocity which may be included to avoid unphysically large velocities.

The Euler-Lagrange equations for \eqref{eq:lagrangian} read
\begin{equation}\label{eq:euler_lagrange}
    \begin{cases}
        \lambda\ddot{\xs} = -\nabla_\xs\beta^2(\Vn(t),\Wm(\xs(t))), \\
        \xs(0) = a, \quad \dot{\xs}(0)=b.
    \end{cases}
\end{equation}
Observe that, if $\lambda=0$, (i.e., no penalization of the velocity is introduced), then \eqref{eq:euler_lagrange} reduces to $\nabla_\xs\beta^2=0$. In this case, the maximization of $\beta^2$ can be achieved via gradient descent techniques. This is the method employed in this work, since large velocities were not observed in the numerical experiments considered. In other contexts, depending on the test problem, it might be necessary to set $\lambda$ equal to a positive value and discretize \eqref{eq:euler_lagrange} directly. Regardless, any numerical scheme employed for the solution of \eqref{eq:euler_lagrange} requires the expression of the gradient of $\beta^2$ with respect to the position of the sensors. We next explain how to compute it.

Let $\xs\in\bR^{md}$ denote the position of the sensors at a fixed time instant, and let us introduce the Gram matrices $\bA(\xs)\in\bR^{m\times m}$ and $\bB(\xs)\in\bR^{m\times n}$, of entries
\begin{equation}\label{eq:B_and_A}
\begin{aligned}
    & (\bA(\xs))_{i,j}:=\bG(\Wm(\xs),\Wm(\xs))_{i,j}=\inprodV{\omega_i(\xs_i)}{\omega_j(\xs_j)}, &\; 1\leq i,j\leq m, \\
    & (\bB(\xs))_{i,s}:=\bG(\Wm(\xs),\Vn)_{i,s} = \inprodV{\omega_i(\xs_i)}{v_s}, &\; 1\leq i\leq m,\, 1\leq s\leq n,
    \end{aligned}
\end{equation}
where $\{v_1,\dots,v_n\}$ is an orthonormal basis of $\Vn$. Then $\beta^2(\Vn,\Wm(\xs))$ is the smallest eigenvalue of the $n\times n$ matrix \cite{Mula2022}
\begin{equation*}
    \bM(\Vn,\Wm(\xs)) := \bB(\xs)^\top\bA(\xs)^{-1}\bB(\xs).
\end{equation*}
For each $\xs$, if $c(\xs)\in\bR^n$ is the unit eigenvector corresponding to the smallest eigenvalue of $\bM(\Vn,\Wm(\xs))$, then
\begin{equation}\label{eq:eigenvalue_prob}
    \bM(\Vn,\Wm(\xs))c(\xs) = \beta^2(\xs)c(\xs),
\end{equation}
which entails
\begin{equation}\label{eq:betasq}
    \beta^2(\xs) = c^\top(\xs)\bM(\Vn,\Wm(\xs))c(\xs).
\end{equation}
Our goal is to obtain an expression for $\nabla_\xs\beta^2(\xs)$. For ease of exposition, we start from the one-dimensional case, that is $d=1$, and consider the general case at the end of the section. Starting from \eqref{eq:betasq} with $\xs=\xs^{(1)}=(\xs_1,\dots,\xs_m)\in\bR^m$, the partial derivative of $\beta^2$ with respect to the $i$th position $\xs_i$ is
\begin{equation}\label{eq:partial_xi_product}
    \frac{\partial\beta^2}{\partial \xs_i}=\frac{\partial c^\top}{\partial \xs_i}\bM(\Vn,\Wm)c+c^\top\frac{\partial\bM(\Vn,\Wm)}{\partial \xs_i}c+c^\top\bM(\Vn,\Wm)\frac{\partial c}{\partial \xs_i},
\end{equation}
where we have dropped all dependencies on $\xs$ for simplicity of notation.
Using \eqref{eq:eigenvalue_prob} yields
\begin{equation*}
    \frac{\partial c^\top}{\partial \xs_i}\bM(\Vn,\Wm)c=\beta^2\frac{\partial c^\top}{\partial \xs_i}c=\beta^2\frac{1}{2}\frac{\partial\lVert c\rVert^2}{\partial \xs_i}=0
\end{equation*}
so that the first and third term in \eqref{eq:partial_xi_product} vanish. We then focus on the partial derivatives of the entries of $\bM$
\begin{equation*}
    c^\top\frac{\partial\bM(\Vn,\Wm)}{\partial \xs_i}c = c^\top\frac{\partial\bB^\top}{\partial \xs_i}\bA^{-1}\bB c + c^\top\bB^\top\frac{\partial\bA^{-1}}{\partial \xs_i}\bB c + c^\top\bB^\top\bA^{-1}\frac{\partial\bB}{\partial \xs_i}c.
\end{equation*}
The first and third terms are equal. In addition, since
$
    \displaystyle\frac{\partial\bA^{-1}}{\partial \xs_i}=-\bA^{-1}\displaystyle\frac{\partial\bA}{\partial \xs_i}\bA^{-1},
$
the second term in \eqref{eq:partial_xi_product} is
\begin{equation}\label{eq:partial_xi}
    c^\top\frac{\partial\bM(\Vn,\Wm)}{\partial \xs_i}c = 2c^\top\frac{\partial \bB^\top}{\partial \xs_i}\bA^{-1}\bB c-c^\top \bB^\top\bA^{-1}\frac{\partial\bA}{\partial \xs_i}\bA^{-1}\bB c.
\end{equation}
Using the fact that, for any $j=1,\dots,m$, $\omega_j$ only depends on the position of the $j$th sensor, the derivative of the $(j,s)$th entry of $\bB$, with $s=1,\dots,n$, reads
\begin{equation*}
    \left(\frac{\partial \bB}{\partial \xs_i}\right)_{j,s}=\inprodV{\frac{\partial\omega_j}{\partial \xs_i}}{v_s}=\inprodV{\frac{\partial\omega_j}{\partial \xs_j}}{v_s}\delta_{ij}
\end{equation*}
where $\delta_{ij}$ is the Kronecker delta. Analogously, differentiating the $(j,k)$th entry of $\bA$, for $j,k=1,\dots,m$, results in
\begin{equation*}
    \left(\frac{\partial\bA}{\partial \xs_i}\right)_{j,k}=\inprodV{\frac{\partial\omega_j}{\partial \xs_i}}{\omega_k}+\inprodV{\omega_j}{\frac{\partial\omega_k}{\partial \xs_i}}=\inprodV{\frac{\partial\omega_j}{\partial \xs_j}}{\omega_k}\delta_{ji}+\inprodV{\omega_j}{\frac{\partial\omega_k}{\partial \xs_k}}\delta_{ik}.
\end{equation*}
Next, we introduce the matrices $\bB_D\in\bR^{m\times n}$ and $\bA_D\in\bR^{m\times m}$ of entries
\begin{equation}\label{eq:BD_and_AD}
    (\bB_D)_{j,s} :=\inprodV{\frac{\partial\omega_j}{\partial \xs_j}}{v_s}\quad\text{ and }\quad (\bA_D)_{j,k}:=\inprodV{\frac{\partial\omega_j}{\partial \xs_j}}{\omega_k},
\end{equation}
for any $j,k=1,\dots,m$, and $s=1,\dots,n$. Then \eqref{eq:partial_xi} can be written as
\begin{align*}
    c^\top\frac{\partial\bM(\Vn,\Wm)}{\partial \xs_i}c&=2(\bB_Dc)_i(\bA^{-1}\bB c)_i-2(\bA^{-1}\bB c)_i(\bA_D\bA^{-1}\bB c)_i=\\&=2[(\bA^{-1}\bB c)\odot(\bB_Dc-\bA_D\bA^{-1}\bB c)]_i,
\end{align*}
where $\odot$ denotes the Hadamard product of vectors. We finally obtain the expression
\begin{equation*}
    \nabla_\xs\beta^2(\xs) = 2(\bA^{-1}\bB c)\odot(\bB_Dc-\bA_D\bA^{-1}\bB c).
\end{equation*}

In the general case $d\geq1$, we have
\begin{equation}\label{eq:gradbetasq}
    \nabla_\xs\beta^2(\xs) = (\nabla_{\xs^{(1)}}\beta^2(\xs), \dots, \nabla_{\xs^{(d)}}\beta^2(\xs))\in\bR^{md},
\end{equation}
and each component $\nabla_{\xs^{(\ell)}}\beta^2(\xs)\in\bR^m$, with $\ell=1,\dots,d$, is obtained by repeating the derivation above for each coordinate $\xs_i^{(\ell)}$ of the $i$th sensor.
This results in
\begin{equation}\label{eq:gradbetasq_xl}
    \nabla_{\xs^{(\ell)}}\beta^2(\xs) = 2(\bA^{-1}\bB c)\odot(\bB^{(\ell)}_Dc-\bA^{(\ell)}_D\bA^{-1}\bB c),
\end{equation}
where $\bA_D^{(\ell)}\in\bR^{m\times m}$ and $\bB_D^{(\ell)}\in\bR^{m\times n}$ are defined as in \eqref{eq:BD_and_AD} with $\xs_j$ replaced by $\xs_j^{(\ell)}$.

We remark that, throughout this section, we implicitly assumed that $\beta^2$ is a smooth function of $\xs$, so that its gradient is well defined. This is always the case when there is no crossing of the eigenvalues of $\bM$ as $\xs$ varies. Moreover, as noted at the beginning of this section, problem \eqref{eq:lagrangian} reduces to the maximization of the smallest eigenvalue of $\bM$ if $\lambda=0$, that is the case considered in this work. This is a non-concave problem in general, and extra care is required to avoid local maxima. An in-depth analysis is out of the scope of this paper and, in practice, no pathological situations were encountered in the numerical tests.

\subsection{Summary of the algorithm and computational complexity}\label{sec:alg_and_complexity}
The procedure for the dynamical placement of the sensors is summarized in \Cref{alg:grad_desc}, for the case $\lambda=0$ in \eqref{eq:lagrangian}. It takes as input the current position of the sensors and the approximation space and returns the updated location of the sensors by solving \eqref{eq:euler_lagrange} with a gradient descent strategy.

\begin{algorithm}[H]
    \caption{Update the position of the sensors}\label{alg:grad_desc}
    \begin{algorithmic}[1]
        \Procedure{$\xs$ = \textsc{Sensors update}}{$\xs$, $\Vn$}
        \State $l\gets 0$
        \While {$l<l_{\text{max}}$}
        \State Compute the Riesz representers $\{\omega_i(\xs_i)\}_{i=1}^m$ spanning $\Wm(\xs)$
        \State Compute $\bB=\bG(\Wm(\xs),\Vn)$, $\bA=\bG(\Wm(\xs),\Wm(\xs))$ as in \eqref{eq:B_and_A}
        \State $\bM\gets\bB^\top{\bA}^{-1}\bB$
        \State Compute the smallest eigenvalue $\beta^2(\xs)$ of $\bM$ and its eigenvector $c$
        \State Compute $\bB_D$ and $\bA_D$ as in \eqref{eq:BD_and_AD}
        \State Compute $\nabla_{\xs}\beta^2(\xs)$ as in \eqref{eq:gradbetasq} and \eqref{eq:gradbetasq_xl}
        \State Determine step size $\alpha$
        \State $\xs\gets\xs+\alpha\nabla_{\xs}\beta^2(\xs)$
        \State $l\gets l+1$
        \EndWhile
        \EndProcedure
    \end{algorithmic}
\end{algorithm}

The computational complexity of \Cref{alg:grad_desc} is
$O(C_S)+O(m^3)+O(mC_R)+O(dm^2C_Q)$ \emph{per iteration}, where:
$C_S$ is the cost of determining the step size $\alpha$ at line 10;
$C_R$ is the cost to numerically solve the boundary value problem to determine the Riesz representers $\omega_i$ for a fixed $i=1,\dots,m$; and $C_Q$ is the number of
operations required to numerically compute the inner products involved in the definition of $\bB$, $\bA$, $\bB_D$ and $\bA_D$. We remark that these costs can be mitigated whenever the expression of the $\{\omega_i\}_{i=1}^m$ is known beforehand, as in the numerical experiments of \Cref{sec:num_exp}.

\section{Application to Hamiltonian dynamics}\label{sec:Ham_dynamics}
In the previous section we have presented a dynamical algorithm to update the position of the sensors, and thus the observation space $\Wm(t)$, assuming that the approximation space $\Vn(t)$ was given at all times $t\in\bR_+$. In this section we consider the opposite perspective and we discuss a strategy to evolve the approximation space $\Vn(t)$. We focus on the case when \eqref{eq:PDE} is a Hamiltonian system. To fix the context, we begin with some definitions.
A Hamiltonian system is a triplet $(M, \cH, \Omega)$ where  $\cH:M\to \bR$ is the so-called Hamiltonian function and $M$ is a manifold endowed with a symplectic structure $\Omega$. In this work we focus on canonical Hamiltonian systems in vector spaces, so we assume that $M=V$ is a Hilbert space. Under this assumption, the pair $(M, \Omega)$ forms a symplectic vector space where $\Omega:V\times V\rightarrow\bR$ is defined as
 $\Omega(u,v):=\inprodV{u}{\cJ(v)}$ for any $u,v\in V$ and $\cJ:V\rightarrow V$ is the linear operator $\cJ(u):=(\up,-\uq)$ for any $u=(\uq,\up)\in V$.
The components $\uq$ and $\up$ of the state variable represent, in Hamiltonian mechanics, the generalized positions in the configuration space and the conjugate momenta, respectively. Hence, $V$ can be written as $V=\Vqp\times\Vqp$ with $\Vqp$ a Hilbert space with inner product $\inprod{\cdot}{\cdot}{\Vqp}$, so that $V$ is naturally endowed with the inner product $\inprodV{u}{v}:=\inprod{u^q}{v^q}{\Vqp}+\inprod{u^p}{v^p}{\Vqp}$ for any $u=(u^q,u^p), v=(v^q,v^p)\in V$.

In this setting, Hamilton's equations are given by
\begin{equation}\label{eq:hamPDE}
    \dot{u} = \cJ \delta \cH(u),
\end{equation}
where $\delta\cH$ is the functional derivative of the Hamiltonian with respect to the state, defined as
$\inprodV{\delta\cH(u)}{v}:= \frac{d}{d\varepsilon}\cH(u+\varepsilon v)\vert_{\varepsilon=0}$ for all $u,v\in V$.
Along solution trajectories, we have conservation of the Hamiltonian; indeed
\begin{equation}\label{eq:hamcons}
    \frac{\rd}{\dt}\cH(u(t)) = \delta \cH(u(t))\dot u(t) = \delta \cH(u)\cJ \delta \cH(u) =0
\end{equation}
due to the skew-symmetry of $\cJ$. The Hamiltonian is often associated with the energy of the system.
Assuming that the Hamiltonian is parameterized by $\theta\in \Theta$, and denoting by $\cH_\theta:V\to\bR$ the corresponding function, we consider the set of solutions $\cM=\cup_{t\geq 0}\cM(t)$ generated by evolutions of the type \eqref{eq:PDE} where $\cP_\theta(u) = \cJ \delta\cH_\theta(u)$.

\subsection{Symplectic dynamical approximation}\label{sec:sdlr}
A good approximation of the set of solutions $\cM(t)$ associated with the Hamiltonian system \eqref{eq:hamPDE} should retain the symplecticity of the approximate trajectories. To this aim we construct approximation spaces that are, for any $t\in\bR_+$, finite-dimensional symplectic vector spaces. Hence, their dimension is even and we denote them by $\Vr(t)$. It is possible to build local-in-time approximations $\Vr(t)$ of the set of solutions $\cM(t)$ associated with the Hamiltonian system \eqref{eq:hamPDE} via an approach based on dynamical low-rank approximation \cite{MN17,P19}. To describe this paradigm, we need to slightly refine the definition of $\cM(t)$ that we gave in \eqref{eq:solution-set-t}. Here, we consider the solution as a function of space and parameter and we introduce the set
\begin{equation*}
    \widehat{\cM}(t) = \{ (\theta, x)\in \Theta\times D \mapsto u(t, x, \theta):\,\mbox{$u$ is solution of \eqref{eq:PDE}} \} \subseteq L^2(\Theta; V).
\end{equation*}
We seek to approximate the elements of $\widehat{\cM}(t)$ on a subspace of the form $\Mr(t) = C_{\Nr}(t)\otimes \Vr(t) \subseteq L^2(\Theta; V)$ where $C_{\Nr}(t)\subset(L^2(\Theta))^{\Nr}$ and $\Vr(t)$ is a linear, time-dependent subspace of $V$ of dimension $2n$. More in detail, we consider the approximation
\begin{equation}\label{eq:udlr}
    u(t, x, \theta) \approx u_{2n}(t, x, \theta) \coloneqq \sum_{i=1}^{2n} c_i(t, \theta) v_i (t, x), \quad \forall\, (t, x, \theta)\in \bR_+\times D\times \Theta,
\end{equation}
where we assume that, for any $t\in\bR_+$, $\bV(t):=(v_1(t), \dots, v_{2n}(t))$ forms a basis of $\Vr(t)$, and we adopt the shorthand notation $v_i(t)=v_i(t,\cdot)$. To ensure preservation of the symplectic structure of the phase space of the approximate dynamics, we need $\Vr(t)$ to be a symplectic subspace of $(V, \Omega)$. For this, we require $\bV(t)$ to be an orthosymplectic basis, that is, $\bV(t)$ belongs to the set
\begin{equation*}
\begin{aligned}
    \Vcal(\Nrh,V):=\left\{\right.(v_1,\ldots,v_{\Nr})\in V^{\Nr}:\, &\,\Omega(v_i,v_j)=(\J{\Nr})_{i,j}\;\mbox{and }\\ & \left.\inprodV{v_i}{v_j}=\delta_{ij},\;\forall\, i,j=1,\ldots,\Nr\right\},
\end{aligned}
\end{equation*}
where $\J{\Nr}$ is the matrix representation of the operator $\Jcal$ in a canonical coordinates system, defined as
\begin{equation*}
    \J{\Nr} \coloneqq
	\begin{pmatrix}
		 \Zrm_{\Nrh} & \Idm_{\Nrh} \\ -\Idm_{\Nrh} & \Zrm_{\Nrh} \\
	\end{pmatrix}\in\R{\Nr}{\Nr},
\end{equation*}
where $\Idm_{\Nrh}, 0_{\Nrh} \in\R{\Nrh}{\Nrh}$ are the identity and zero matrices, respectively. Note that $\Vcal(\Nrh,V)$ is a subset of the Stiefel manifold of $\Nr$-dimensional symplectic bases of the symplectic vector space $(V,\Omega)$.

Moreover, for a given $t\geq0$, we require the coefficients in \eqref{eq:udlr} to belong to the space
\begin{equation*}
    \Ccal(\Nr,L^2(\Theta)):=\left\{\bC=(c_1,\dots, c_{2n})\in (L^2(\Theta))^{\Nr}\cond \rank(\bS(\bC)) =\Nr\right\},
\end{equation*}
where
\begin{equation*}
    \bS(\bC):=\int_{\Theta}\bC^{\top}\bC\, d\prm +\J{\Nr}^{\top}\left( \int_{\Theta}\bC^{\top}\bC\, d\prm\right)\,\J{\Nr}\in\bR^{\Nr\times\Nr} \quad \forall\, \bC \in (L^2(\Theta))^{2n}.
\end{equation*}
This construction implies that, for every $t\in\bR_+$, the elements of $\widehat{\cM}(t) \subset L^2(\Theta; V)$ are approximated as
\begin{equation}\label{eq:fact}
    u_{2n}(t) = \bV(t)\bC^{\top}(t)\in \Mr(t),
\end{equation}
where
\begin{equation*}
    \Mr(t)\coloneqq\left\{\bV(t)\bC^{\top}(t)\cond \bV(t)\in\Vcal(\Nrh,V),\;\bC(t)\in\Ccal(\Nr,L^2(\Theta))\right\}\subset (L^2(\Theta))^{\Nr}\otimes V^{\Nr}.
\end{equation*}
The fact that we require $\bC\in\Ccal(\Nr,L^2(\Theta))$ ensures uniqueness of the decomposition in the sense that, if $\bV_1\bC^{\top}=\bV_2 \bC^{\top}$ for $\bV_1, \bV_2\in \Vcal(\Nrh,V)$, then $\bV_1=\bV_2$.

The space $\Mr(t)$ is a submanifold of $(L^2(\Theta))^{\Nr}\otimes (V, \Omega)^{\Nr}$ and its tangent space at each point $u_{2n}\in \Mr(t)$ can be characterized as
\begin{equation*}
\begin{aligned}
    T_{\ur(t)}\Mr(t):=\left\{\right.& \delta\bV(t)\,\bC^{\top}(t) +\bV(t)\,\delta\bC^{\top}(t)\in (L^2(\Theta)\otimes V)^{\Nr}:\\ &    \delta\bC(t)\in(L^2(\Theta))^{\Nr},\, \delta\bV(t)=(\delta v_1,\ldots,\delta v_{\Nr})\in V^{\Nr}\\ &\left.        \mbox{with}\; \delta\bV(t)\J{\Nr}=\Jcal\delta\bV(t),\,\inprodV{\delta v_i}{v_j}=0\right\}.
\end{aligned}
\end{equation*}

The approximate trajectory $t\in\bR_+\mapsto \ur(t)\in\Mr(t)$ is determined via the Dirac--Frenkel variational principle \cite{Fren34,ML64}, which imposes that the vector field $\dur\in T_{\ur(t)}{\Mr(t)}$ satisfies, at every time $t$,
\begin{equation}\label{eq:DF}
    \int_{\Theta}\Omega(\dur-\cP_{\prm}(\ur),\xi) \,d\theta=0\qquad\forall\, \xi\in T_{\ur(t)}{\Mr(t)}.
\end{equation}
This amounts to projecting the vector field $\cP_{\prm}$ at $\ur(t)$ to the tangent space of $\Mr(t)$ at $\ur(t)$ via a symplectic projection.

Applying the Dirac--Frenkel variational principle \eqref{eq:DF} to problem \eqref{eq:hamPDE} and using the factorization \eqref{eq:fact}, the dynamics of the approximate state $\ur(t)=\bV(t)\bC^{\top}(t)$ on $\Mr$ can be written in terms of evolution equations for the basis $\bV$ and for the coefficients $\bC$. Using the explicit expression for the projection given in \cite{MN17,P19} yields
\begin{equation}\label{eq:DLR}
\left\{
\begin{aligned}
    & \dot{\bV}\,\bS(\bC)  = P_{\Vr^{\perp}(t)}\Big(\int_\Theta(\cP_{\prm}(\ur)\bC+\Jcal\cP_{\prm}(\ur)\bC\J{\Nr}^{\top})\,d\theta\Big), \\
    & \dot{\bC}  = \inprod{P_{\Vr(t)}\cP_{\prm}(\ur)}{\bV}{(L^2(\Sdom))^{2}},
\end{aligned}\right.
\end{equation}
where $P_{\Vr}$ and $P_{\Vr^{\perp}}$ denote the $V$-orthogonal projection onto $\Vr$ and $\Vr^{\perp}$, respectively, and the initial condition of problem \eqref{eq:DLR} is obtained from the projection of $u_0$ onto $\Mr$.

We conclude this section with a remark on our assumption that $\widehat{\cM}(t) \subset L^2(\Theta; V)$. This hypothesis could be understood both from an analytical and a probabilistic angle. From the analytical perspective, we are simply assuming $L^2$ summability in $\Theta$ of $\widehat{\cM}(t)$ with respect to the Lebesgue measure. Alternatively, one can take a more probabilistic point of view and interpret the assumption as the parameters following a uniform distribution over $\Theta$. This connects with Bayesian approaches, and would motivate to update the underlying distribution over time. In this work, we are rather taking the analytical point of view. As a consequence, integrals over $\Theta$ in \eqref{eq:DLR} are approximated using deterministic quadrature rules, and not sampling techniques. We refer to \Cref{sec:num_exp} for more details on this aspect. The probability-driven perspective will be explored in future works.

\subsection{Properties of the reconstructed flow}\label{sec:Ham_pres}
The construction and evolution of the approximation space $\Vr(t)$ according to \Cref{sec:sdlr} ensures that the reconstructed trajectory $t\mapsto v_t^*(\omega(t))$ obtained in \eqref{eq:vtstar} belongs to $\cup_{t\geq 0}\Vr(t)$, which is the union of symplectic vector spaces, and therefore the symplecticity of the flow is preserved. This is a key aspect in preventing the reconstructed dynamics from developing unphysical oscillations and unstable behaviors. On the other hand, the Hamiltonian associated with the true parameter is not guaranteed to be preserved by the reconstructed flow. Nevertheless, the error in the conservation of the Hamiltonian can be related to the reconstruction error \eqref{eq:Evtstar}. Owing to the definition of $\cM(t)$ \eqref{eq:solution-set-t}, we shall hereon denote this quantity as
\begin{equation*}
    E(v_t^*,\cM(t)) = \Ecal_\Theta(t):=\sup_{\theta\in\Theta}\lVert u_\theta(t)-v_\omega^*(t)\rVert,
\end{equation*}
with the shorthand notation $u_\theta(t):=u(\theta)(t,\cdot)$ and $v^*_\omega(t):=v_t^*(\omega(t))$. If we assume that, for all $\theta\in\Theta$, the Hamiltonian $\Hcal_\theta:V\rightarrow \bR$ is Lipschitz continuous, uniformly in $t$, with Lipschitz constant $L_\theta$, then
\begin{equation*}
    \Ecal_{\Hcal,\Theta}(t):=\sup_{\theta\in\Theta}\lvert\Hcal_\theta(u_\theta(t))-\Hcal_\theta(v^*_\omega(t))\rvert\leq L\Ecal_\Theta(t)\leq L\frac{\text{dist}(\cM(t),\Vn(t))}{\beta(\Vn(t),\Wm(t))}
\end{equation*}
where $L:=\sup_{\theta\in\Theta}L_\theta$ and the second inequality follows from \eqref{eq:Evtstar}. Moreover, if we consider the conservation of the Hamiltonian evaluated on the reconstructed trajectories, we have that
\begin{equation*}
    \Delta\Hcal_\Theta(t) := \sup_{\theta\in\Theta}\lvert\Hcal_\theta(v_\omega^*(t))-\Hcal_\theta(v_\omega^*(0))\rvert\leq L(\Ecal_\Theta(t)+\Ecal_\Theta(0)).
\end{equation*}
This follows by simple triangle inequality
\begin{align*}
    \lvert\Hcal_\theta(v_\omega^*(t))-\Hcal_\theta(v_\omega^*(0))\rvert \leq &\, \lvert\Hcal_\theta(v_\omega^*(t))-\Hcal_\theta(u_\theta(t))\rvert + \lvert\Hcal_\theta(u_\theta(0))-\Hcal_\theta(v_\omega^*(0))\rvert\\ & + \lvert\Hcal_\theta(u_\theta(t))-\Hcal_\theta(u_\theta(0))\rvert,\qquad \forall\,\theta\in\Theta,
\end{align*}
and the fact that the Hamiltonian is conserved along the exact trajectory, see \eqref{eq:hamcons}.
This result states that the Hamiltonian error in the PBDW approximation is bounded by the reconstruction error with a factor depending on the regularity of $\Hcal_\theta$. The hypothesis of $\Hcal_\theta$ being Lipschitz is typically not restrictive, since it follows from the Lipschitz continuity of $\cP_\theta$, which in turn guarantees well-posedness of \eqref{eq:hamPDE}.

\subsection{Summary of the Dyn-PBDW algorithm and computational complexity}
The dynamical state estimation procedure is summarized in \Cref{alg:dse}.

\begin{algorithm}
    \caption{Dynamical state estimation of Hamiltonian systems}\label{alg:dse}
    \begin{algorithmic}[1]
        \Procedure{\,Dyn-PBDW}{$\Vtwon(0)$, $\xs(-\Delta t)$, $T$, $\Delta t$}
        \State $t_0\gets-\Delta t$, $t_1\gets0$
        \While {$t_1\leq T$}
        \State $\xs(t_1) =$ \textsc{Sensors update}$(\xs(t_0),\Vtwon(t_1))$ as in \Cref{alg:grad_desc}
        \State Obtain measurements $\omega(\xs(t_1))$
        \State Compute the PBDW reconstruction as in \Cref{sec:DPBDW-method}
        \State $t_0\gets t_1$, $t_1\gets t_0+\Delta t$
        \State Update the approximation space $\Vtwon(t_0)$ to $\Vtwon(t_1)$ by solving \eqref{eq:DLR}
        \EndWhile
        \EndProcedure
    \end{algorithmic}
\end{algorithm}

The inputs are the approximation space at the initial time $\Vtwon(0)$ and an initial guess for the position of the sensors, denoted by $\xs(-\Delta t)$. At each time step during the simulation, the positions of the sensors are updated at line 4 as described in \Cref{alg:grad_desc}. The expression of $\nabla_{\xs}\beta^2$ in the case of Hamiltonian systems is derived in \Cref{app:gradbeta_Ham}. Next, the reconstructed solution is obtained with the PBDW method described in \Cref{sec:DPBDW-method} starting from measurements at the new positions. The total complexity of the PBDW algorithm is negligible compared to the cost of moving the sensors. Finally, the approximation space is evolved at line 8. To this end, \eqref{eq:DLR} is discretized in space and the parameter space is approximated by the finite set $\Thtrain\subset\Theta$. Assuming that $2N$ is the number of degrees of freedom of the spatial discretization, we obtain semi-discrete systems of size $2N\times2n$ and $2n\times\lvert\Thtrain\rvert$ that are then numerically integrated in time. In this work, we focus on the family of partitioned Runge--Kutta method developed in \cite{P19,HPRi22}, whose complexity is $O(Nn^2s)+C_{\cP}$, where $s$ is the number of stages of the Runge--Kutta scheme and $C_{\cP}$ is the cost of evaluating the nonlinear vector field $\cP_\theta$ for all $\theta\in\Thtrain$. The complexity of this operation depends on both $N$ and $\lvert\Thtrain\rvert$, and it can be mitigated via hyper-reduction techniques \cite{PV22,PV23}. Therefore, owing to the analysis carried out in \Cref{sec:alg_and_complexity}, the total computational complexity of Dyn-PBDW is $O(Nn^2s)+C_{\cP}+l_{\max}\big(O(mC_R)+O(m^3)+O(dm^2C_Q)+C_S\big)$.

\section{Numerical results}\label{sec:num_exp}
In this section we assess the performance of Dyn-PBDW on the (cubic) nonlinear Schr\"odinger equation and the shallow water equations. The computational domain $\Sdom$ is the interval $[-L_x,L_x]$, if $d=1$, and the rectangle $[-L_x,L_x]\times[-L_y,L_y]$, if $d=2$. A high-fidelity approximation of the problem is obtained by means of a second order finite difference discretization on a computational grid of $N_x$ and $N_y$ equispaced points in each direction, and it is solved until the final time $T$ using $N_t$ uniform time steps. The implicit midpoint rule is chosen as the time integration scheme, and the nonlinear system arising at each time step is solved with the Newton method with stopping tolerance $\tau=10^{-10}$. The numerical solution of the high-fidelity system is then regarded as the ground truth for state estimation.

The evolution of the approximation space $\Vr(t)$ is performed as in \Cref{sec:sdlr}. The same spatial discretization as the high-fidelity system is applied to \eqref{eq:DLR}, and $\Theta$ is discretized by a finite set $\Thtrain$ of uniformly selected test parameters. The test parameters in $\Thtrain$ are used to approximate the integrals over $\Theta$ appearing in \eqref{eq:DLR} using a deterministic quadrature rule.
The temporal discretization of the approximate dynamics is done via partitioned Runge--Kutta methods as in \cite{P19}, \cite[Section 4]{HPRi22}.

The setting for all numerical tests is the Hilbert space $V=\Vqp\times\Vqp$ where $\Vqp=L^2(D)$. For $u=(u^q,u^p)\in V$, $2m$ measurements are given in the form of local averages
\begin{equation*}
    \ell_i(u)=
    \begin{cases}
        \ell^q_i(u^q) = \displaystyle\frac{1}{(2\pi\sigma^2)^{\frac{d}{2}}}\int_\Sdom\exp{\left(-\displaystyle\frac{\lVert x-\xs_i\rVert^2}{2\sigma^2}\right)}u^q(x)\,dx & i=1,\dots,m, \\
        \ell^p_i(u^p) = \displaystyle\frac{1}{(2\pi\sigma^2)^{\frac{d}{2}}}\int_\Sdom\exp{\left(-\displaystyle\frac{\lVert x-\xs_{i-m}\rVert^2}{2\sigma^2}\right)}u^p(x)\,dx & i=m+1,\dots,2m,
    \end{cases}
\end{equation*}
where $\xs_i\in D$ is the position of the $i$th sensor. In this setting, the Riesz representer of $\ell^q_i=\ell^p_{i+m}$ is simply
\begin{equation*}
    \omega^q_i(x)=\omega^p_{i+m}(x)=\frac{1}{(2\pi\sigma^2)^{\frac{d}{2}}}\exp{\left(-\frac{\lVert x-\xs_i\rVert^2}{2\sigma^2}\right)} \qquad i=1,\dots,m.
\end{equation*}
In the numerical tests we set $\sigma=0.1$. Qualitatively similar results were observed for other choices of $\sigma$.

In order to assess the quality of the reconstruction, we first define the state estimation error associated to the parameter $\theta\in\Theta$ in the norm of $V$
\begin{equation*}
    \Ese(t,\theta)\coloneqq \lVert u(\theta)(t, \cdot)-v_t^*(\omega(t))\rVert,
\end{equation*}
where $u$ is the solution of the high-fidelity system and $\omega$ is its projection onto $\Wtwom(t)$. We also define the quantities
\begin{equation*}
    \Eproj(t,\theta)\coloneqq \lVert u(\theta)(t, \cdot)-P_{\Vtwon(t)}u(\theta)(t, \cdot)\rVert \quad \text{and} \quad \Ebound(t,\theta)\coloneqq \Eproj(t,\theta)/\beta(t)
\end{equation*}
which represent the projection error of the high-fidelity solution onto the approximation space and the error bound for the PBDW reconstruction, respectively. Since the reconstruction belongs to the approximation space, these errors satisfy $\Eproj(t,\theta)\leq\Ese(t,\theta)\leq\Ebound(t,\theta)$ for any parameter $\theta$. In the numerical experiments, state estimation is performed on a finite subset $\Thtest\subset \Theta$ of cardinality $\lvert\Thtest\rvert$. To study reconstruction errors, we compute
\begin{equation*}
    \Esets(t) \coloneqq  \max_{\theta\in\Thtest}\Ese(t,\theta), \qquad \Eprojts(t)\coloneqq \max_{\theta\in\Thtest}\Eproj(t,\theta), \qquad \Eboundts(t)\coloneqq \max_{\theta\in\Thtest}\Ebound(t,\theta),
\end{equation*}
and the integrals involved in $V$-norm computations are approximated by means of a suitable quadrature rule. These quantities are estimates of the errors introduced in Section \ref{sec:DPBDW}: $\Esets$ is an estimate of the reconstruction error $E(v_t^*, \cM(t))$ defined in \eqref{eq:Evtstar}; $\Eprojts$ approximates $\dist(\cM(t), \Vn(t))$ which is a lower bound of $E(v_t^*, \cM(t))$; and $\Eboundts$ estimates the upper bound in \eqref{eq:Evtstar}. In particular, we have $\Eprojts\leq \Esets\leq \Eboundts$. When needed, we shall distinguish between the errors obtained with the static and dynamic placement of the sensors by denoting them as $\Esetsstat$ and $\Esetsdyn$, respectively.

In this work we are also interested in quantifying the difference between the Hamiltonian evaluated at the ground truth and at the reconstructed solution. Analogously to \Cref{sec:Ham_pres} we introduce the Hamiltonian error
\begin{equation}\label{eq:eHam}
    \eHts(t) \coloneqq  \max_{\theta\in\Thtest}\lvert\Hcal_{\theta}(u(\theta)(t, \cdot))-\Hcal_{\theta}(v_t^*(\omega(t)))\rvert,
\end{equation}
and the error in the conservation of the Hamiltonian
\begin{equation}\label{eq:dHam}
    \dHts(t) \coloneqq  \max_{\theta\in\Thtest}\lvert\Hcal_{\theta}(\mathfrak{u}(t))-\Hcal_{\theta}(\mathfrak{u}(0))\rvert,
\end{equation}
where $\mathfrak{u}(t)$ is either the high-fidelity solution $u(\theta)(t,\cdot)$ or the reconstruction $v_t^*(\omega(t))$.

The MATLAB\textsuperscript{\textregistered} code used to obtain the results presented in this section is available at \url{https://github.com/fvismara/Dyn-PBDW/tree/main}.

\subsection{1D Schr\"odinger equation}
We first consider the one-dimensional, cubic Schr\"odinger equation
\begin{equation}\label{eq:1dNLS}
    \begin{cases}
        \imath \psi_t+\psi_{xx}+\epsilon\lvert \psi\rvert^2\psi=0
         &\quad x\in\Sdom, \; t\in(0,T], \\
        \psi(x,0) = \displaystyle\frac{\sqrt{2}}{\cosh{(\alpha x)}}\exp\left(\imath\displaystyle\frac{x}{2}\right) &\quad x\in\Sdom,
    \end{cases}
\end{equation}
with periodic boundary conditions. The subscripts $t$ and $x$ denote partial differentiation with respect to the corresponding variable. The solution depends on the parameter vector $\theta=(\alpha,\epsilon)\in\Theta=[0.98,1.1]^2\subset\bR^2$. A Hamiltonian formulation is obtained by writing the complex-valued wave function $\psi$ in terms of its real and imaginary parts as $\psi(x,t)=q(x,t)+\imath p(x,t)$ and recasting \eqref{eq:1dNLS} in terms of evolution equations for $q$ and $p$:
\begin{equation*}
    \begin{cases}
        q_t = -p_{xx}-\epsilon(q^2+p^2)p, \\
        p_t = q_{xx}+\epsilon(q^2+p^2)q.
    \end{cases}
\end{equation*}
This is a Hamiltonian system with Hamiltonian
\begin{equation*}
    \Hcal_\theta(q,p)=\frac{1}{2}\int_\Sdom \left[q_x^2+p_x^2-\frac{\epsilon}{2}(q^2+p^2)^2\right]\,dx. 
\end{equation*}

\begin{figure}[ht]
    \centering
    \includegraphics[width=0.95\textwidth]{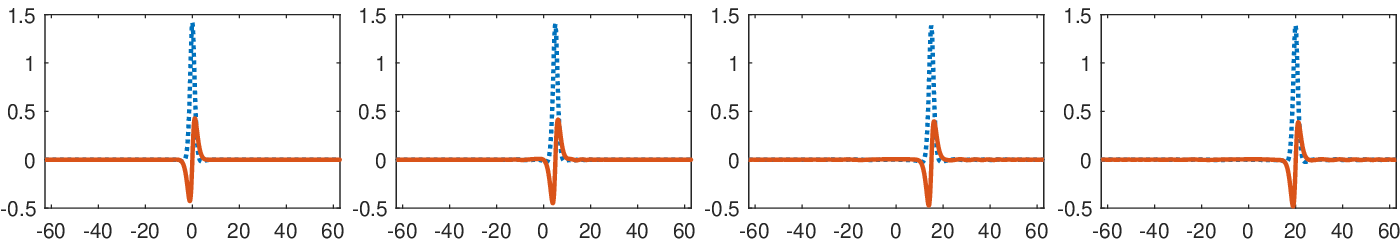}
    \caption{\footnotesize 1D Schr\"odinger equation. Real part $q(x,t)$ (dotted) and imaginary part $p(x,t)$ (solid) of the high-fidelity solution at times $t=0$, $t=5$, $t=15$ and $t=20$ for $\theta=(1.04,1.04)$.}\label{fig:1DNLS_fullsol}
\end{figure}

In the numerical experiments we set $L_x=20\pi$, $N_x=1000$, $T=20$ and $N_t=20000$. We show in \Cref{fig:1DNLS_fullsol} the real and imaginary part of the high-fidelity solution corresponding to the parameter $\theta=(1.04,1.04)$ at four time instants. The solution is characterized by a main hump moving to the right, with small oscillations originating around it. A qualitatively similar behaviour is observed for all test parameters.

The approximation space $\Vr$ in this test has dimension $2n=8$, and its evolution is computed using $\lvert\Thtrain\rvert=10^2=100$ uniform parameters. Then, we perform state estimation for $\lvert\Thtest\rvert=9^2=81$ test parameters.

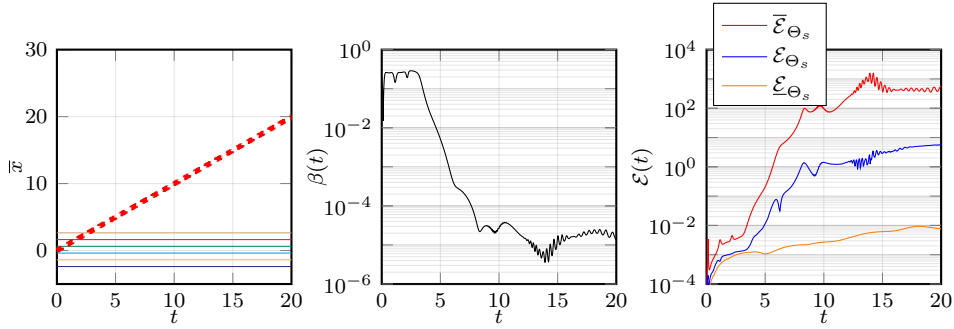
\begin{figure}[ht]
\centering
\begin{tikzpicture}
    \begin{groupplot}[
      group style={group size=3 by 1, horizontal sep=1.2cm},
      width=4.7cm, height=4.7cm
    ]
    \nextgroupplot[xlabel={$t$},
                  ylabel={$\xs$},
                  xlabel style = {yshift=.2cm},
                  ylabel style = {yshift=-.2cm},
                  axis line style = thick,
                  grid=both,
                  grid style = {gray,opacity=0.2},
                  xmin=0, xmax=20,
                  ymin=-5, ymax=30,
                  xlabel style={font=\footnotesize},
                  ylabel style={font=\footnotesize},
                  x tick label style={font=\footnotesize},
                  y tick label style={font=\footnotesize}]
        \addplot+[color=Blue,mark=none] table[x=t,y=s1x] {1DNLS_noadapt_pointwise.txt};
        \addplot+[color=Dandelion,mark=none] table[x=t,y=s2x] {1DNLS_noadapt_pointwise.txt};
        \addplot+[color=Cerulean,mark=none] table[x=t,y=s3x] {1DNLS_noadapt_pointwise.txt};
        \addplot+[color=PineGreen,mark=none] table[x=t,y=s4x] {1DNLS_noadapt_pointwise.txt};
        \addplot+[color=OrangeRed,mark=none] table[x=t,y=s5x] {1DNLS_noadapt_pointwise.txt};
        \addplot+[color=Tan,solid,mark=none] table[x=t,y=s6x] {1DNLS_noadapt_pointwise.txt};
        \addplot+[color=red,dashed,mark=none,line width=1.5pt] table[x=t,y=bumpx] {1DNLS_noadapt_pointwise.txt};
    \nextgroupplot[xlabel={$t$},
                  ylabel={$\beta(t)$},
                  xlabel style = {yshift=.2cm},
                  ylabel style = {yshift=-.3cm},
                  axis line style = thick,
                  grid=both,
                  ymode=log,
                  grid style = {gray,opacity=0.2},
                  xmin=0, xmax=20,
                  ymin=1e-06, ymax=1,
                  ytick={1e-06,1e-05,1e-04,1e-03,1e-02,1e-01,1e+00},
                  yticklabels={},
                  extra y ticks={1e-06,1e-04,1e-02,1e+00},
                  yminorticks=false,
                  xlabel style={font=\footnotesize},
                  ylabel style={font=\footnotesize},
                  x tick label style={font=\footnotesize},
                  y tick label style={font=\footnotesize}]
        \addplot+[color=black,mark=none] table[x=t,y=beta] {1DNLS_noadapt_pointwise.txt};
    \nextgroupplot[xlabel={$t$},
                  ylabel={$\Ecal(t)$},
                  xlabel style = {yshift=.2cm},
                  ylabel style = {yshift=-.3cm},
                  axis line style = thick,
                  ytick={1e-04,1e-03,1e-02,1e-01,1e+00,1e+01,1e+02,1e+03,1e+04},
                  yticklabels={},
                  extra y ticks = {1e-04,1e-02,1e+00,1e+02,1e+04},
                  yminorticks=false,
                  grid=both,
                  ymode=log,
                  grid style = {gray,opacity=0.2},
                  xmin=0, xmax=20,
                  ymin=1e-04, ymax=1e+04,
                  xlabel style={font=\footnotesize},
                  ylabel style={font=\footnotesize},
                  x tick label style={font=\footnotesize},
                  y tick label style={font=\footnotesize},
                  legend style={font=\footnotesize},
                  legend cell align={left},
                  legend columns = 1,
                  legend style={at={(0.27,1.2)},anchor=north}]
        \addplot+[color=red,mark=none] table[x=t,y=eb] {1DNLS_noadapt_pointwise.txt};
        \addplot+[color=blue,mark=none] table[x=t,y=ese] {1DNLS_noadapt_pointwise.txt};
        \addplot+[color=orange,mark=none] table[x=t,y=ep] {1DNLS_noadapt_pointwise.txt};
        \legend{{$\Eboundts$},{$\Esets$},{$\Eprojts$}};
    \end{groupplot}
\end{tikzpicture}
\caption{\footnotesize 1D Schr\"odinger equation. Static local averages of width $\sigma=0.1$, $m=6$. Left: $x$ coordinate of the sensors over time; the red dashed line denotes the position of the main hump of the high-fidelity solution $\lvert\psi(t,\theta)\rvert$ for $\theta=(1.04,1.04)$. Center: evolution of the inf-sup constant $\beta$. Right: evolution of the numerical errors with $\lvert\Thtest\rvert=81$ test parameters.}\label{fig:1DNLS_noadapt_pointwise}
\end{figure}

With the first test we aim at showing the need for a dynamic selection of the sensor positions. We consider a small number of sensors, $m=6$, which are initially placed around the main hump and do not move over time, see \Cref{fig:1DNLS_noadapt_pointwise}, left. Here the sensors do not follow the evolution of the numerical solution, and, after some time, they only measure the small oscillations around the hump. Therefore, while being large at the initial stages, the value of $\beta$ decreases as time evolves, its value at the final time being around $10^{-5}$. The effect is that the reconstruction error deteriorates over time, as shown in \Cref{fig:1DNLS_noadapt_pointwise} on the right.

\begin{figure}[ht]
\centering
\begin{tikzpicture}
    \begin{groupplot}[
      group style={group size=3 by 1, horizontal sep=1.2cm},
      width=4.7cm, height=4.7cm
    ]
    \nextgroupplot[xlabel={$t$},
                  ylabel={$\xs$},
                  xlabel style = {yshift=.2cm},
                  ylabel style = {yshift=-.2cm},
                  axis line style = thick,
                  grid=both,
                  grid style = {gray,opacity=0.2},
                  xmin=0, xmax=20,
                  ymin=-5, ymax=30,
                  xlabel style={font=\footnotesize},
                  ylabel style={font=\footnotesize},
                  x tick label style={font=\footnotesize},
                  y tick label style={font=\footnotesize}]
        \addplot+[color=Blue,mark=none] table[x=t,y=s1x] {1DNLS_adapt_pointwise.txt};
        \addplot+[color=Dandelion,mark=none] table[x=t,y=s2x] {1DNLS_adapt_pointwise.txt};
        \addplot+[color=Cerulean,mark=none] table[x=t,y=s3x] {1DNLS_adapt_pointwise.txt};
        \addplot+[color=PineGreen,mark=none] table[x=t,y=s4x] {1DNLS_adapt_pointwise.txt};
        \addplot+[color=OrangeRed,mark=none] table[x=t,y=s5x] {1DNLS_adapt_pointwise.txt};
        \addplot+[color=Tan,solid,mark=none] table[x=t,y=s6x] {1DNLS_adapt_pointwise.txt};
        \addplot+[color=red,dashed,mark=none,line width=1.5pt] table[x=t,y=bumpx] {1DNLS_adapt_pointwise.txt};
    \nextgroupplot[xlabel={$t$},
                  ylabel={$\beta(t)$},
                  xlabel style = {yshift=.2cm},
                  ylabel style = {yshift=-.3cm},
                  axis line style = thick,
                  grid=both,
                  ymode=log,
                  grid style = {gray,opacity=0.2},
                  xmin=0, xmax=20,
                  ymin=1e-06, ymax=1,
                  ytick={1e-06,1e-05,1e-04,1e-03,1e-02,1e-01,1e+00},
                  yticklabels={},
                  extra y ticks = {1e-06,1e-04,1e-02,1e+00},
                  yminorticks = false,
                  xlabel style={font=\footnotesize},
                  ylabel style={font=\footnotesize},
                  x tick label style={font=\footnotesize},
                  y tick label style={font=\footnotesize}]
        \addplot+[color=black,mark=none] table[x=t,y=beta] {1DNLS_adapt_pointwise.txt};
    \nextgroupplot[xlabel={$t$},
                  ylabel={$\Ecal(t)$},
                  xlabel style = {yshift=.2cm},
                  ylabel style = {yshift=-.3cm},
                  axis line style = thick,
                  ytick = {1e-03,1,1e+03},
                  grid=both,
                  ymode=log,
                  grid style = {gray,opacity=0.2},
                  xmin=0, xmax=20,
                  ymin=1e-04, ymax=1e+04,
                  ytick={1e-04,1e-03,1e-02,1e-01,1e+00,1e+01,1e+02,1e+03,1e+04},
                  yticklabels = {},
                  extra y ticks = {1e-04,1e-02,1e+00,1e+02,1e+04},
                  yminorticks = false,
                  xlabel style={font=\footnotesize},
                  ylabel style={font=\footnotesize},
                  x tick label style={font=\footnotesize},
                  y tick label style={font=\footnotesize},
                  legend style={font=\footnotesize},
                  legend cell align={left},
                  legend columns = 1,
                  legend style={at={(0.27,0.98)},anchor=north}]
        \addplot+[color=red,mark=none] table[x=t,y=eb] {1DNLS_adapt_pointwise.txt};
        \addplot+[color=blue,mark=none] table[x=t,y=ese] {1DNLS_adapt_pointwise.txt};
        \addplot+[color=orange,mark=none] table[x=t,y=ep] {1DNLS_adapt_pointwise.txt};
        \legend{{$\Eboundts$},{$\Esets$},{$\Eprojts$}};
    \end{groupplot}
\end{tikzpicture}
\caption{\footnotesize 1D Schr\"odinger equation. Dynamic local averages of width $\sigma=0.1$, $m=6$. Left: $x$ coordinate of the sensors over time; the red dashed line denotes the position of the main hump of the high-fidelity solution $\lvert\psi(t,\theta)\rvert$ for $\theta=(1.04,1.04)$. Center: evolution of the inf-sup constant $\beta$. Right: evolution of the numerical errors with $\lvert\Thtest\rvert=81$ test parameters.}\label{fig:1DNLS_adapt_pointwise}
\end{figure}
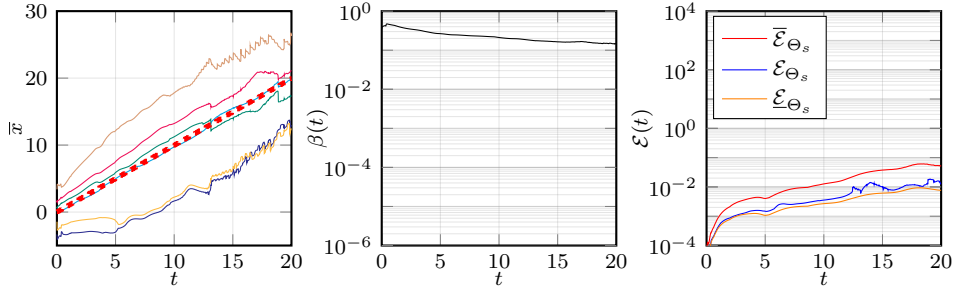

\Cref{fig:1DNLS_adapt_pointwise} shows the outcome of Dyn-PBDW: the locations of the sensors are updated in time, and they closely track the support of the solution. Thanks to this, the value of $\beta$ does not exhibit a significant decay and the reconstruction error $\cE_{\Thtest}$ remains close to the best approximation error $\underline{\cE}_{\Thtest}$ at all times.

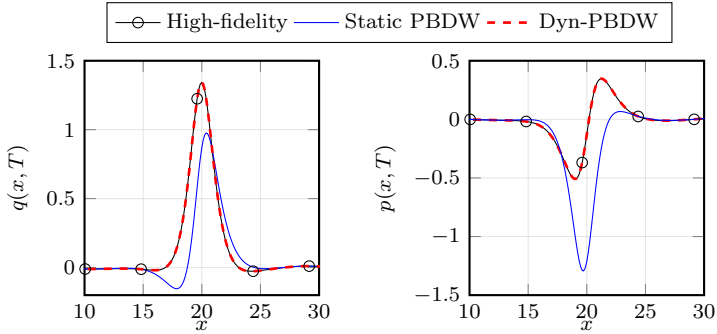
\begin{figure}[ht]
\centering
\begin{tikzpicture}
    \begin{groupplot}[
      group style={group size=2 by 1, horizontal sep=2cm},
      width=4.7cm, height=4.7cm
    ]
    \nextgroupplot[xlabel={$x$},
                  ylabel={$q(x,T)$},
                  xlabel style = {yshift=.2cm},
                  axis line style = thick,
                  grid=both,
                  grid style = {gray,opacity=0.2},
                  xmin=10, xmax=30,
                  ymin=-0.2, ymax=1.5,
                  xlabel style={font=\footnotesize},
                  ylabel style={font=\footnotesize},
                  x tick label style={font=\footnotesize},
                  y tick label style={font=\footnotesize},
                  legend style={font=\footnotesize},
                  legend cell align={left},
                  legend columns = 3,
                  legend style={at={(1.3,1.25)},anchor=north}]
        \addplot+[color=black,mark=o,mark repeat={38}] table[x=xgrid,y=qex] {1DNLS_recsol_a_vs_na.txt};
        \addplot+[color=blue,mark=none] table[x=xgrid,y=qna] {1DNLS_recsol_a_vs_na.txt};
        \addplot+[color=red,dashed,mark=none,line width=1] table[x=xgrid,y=qa] {1DNLS_recsol_a_vs_na.txt};
        \legend{{High-fidelity},{Static PBDW},{Dyn-PBDW}};
    \nextgroupplot[xlabel={$x$},
                  ylabel={$p(x,T)$},
                  xlabel style = {yshift=.2cm},
                  axis line style = thick,
                  grid=both,
                  grid style = {gray,opacity=0.2},
                  xmin=10, xmax=30,
                  ymin=-1.5, ymax=0.5,
                  xlabel style={font=\footnotesize},
                  ylabel style={font=\footnotesize},
                  x tick label style={font=\footnotesize},
                  y tick label style={font=\footnotesize}]
        \addplot+[color=black,mark=o,mark repeat={38}] table[x=xgrid,y=pex] {1DNLS_recsol_a_vs_na.txt};
        \addplot+[color=red,dashed,mark=none,line width=1] table[x=xgrid,y=pa] {1DNLS_recsol_a_vs_na.txt};
        \addplot+[color=blue,mark=none] table[x=xgrid,y=pna] {1DNLS_recsol_a_vs_na.txt};
    \end{groupplot}
\end{tikzpicture}
\caption{\footnotesize 1D Schr\"odinger equation. Numerical solution of the high-fidelity model corresponding to $\theta=(1.0933,1.0933)$ and comparison with the statically and dynamically reconstructed solutions at the final time. For the sake of visualization, only the interval $[10,30]\subset\Sdom$ is shown.}\label{fig:1DNLS_recsol}
\end{figure}

A comparison between the reconstructed solutions obtained with the static and dynamic placement of the sensors is presented for one test parameter in \Cref{fig:1DNLS_recsol}. The dynamically reconstructed solution is visually indistinguishable from the ground truth, while the reconstruction obtained with static sensors fails to correctly capture the location and amplitude of the main profiles.

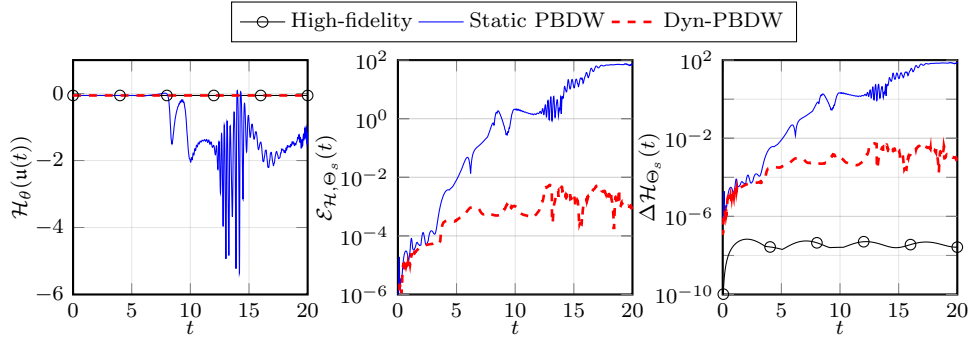
\begin{figure}[ht]
\centering
\begin{tikzpicture}
    \begin{groupplot}[
      group style={group size=3 by 1, horizontal sep=1.2cm},
      width=4.7cm, height=4.7cm
    ]
    \nextgroupplot[xlabel={$t$},
                  ylabel={$\Hcal_\theta(\mathfrak{u}(t))$},
                  xlabel style = {yshift=.2cm},
                  ylabel style = {yshift=-.2cm},
                  axis line style = thick,
                  grid=both,
                  grid style = {gray,opacity=0.2},
                  xmin=0, xmax=20,
                  ymin=-6, ymax=1,
                  xlabel style={font=\footnotesize},
                  ylabel style={font=\footnotesize},
                  x tick label style={font=\footnotesize},
                  y tick label style={font=\footnotesize},
                  legend style={font=\footnotesize},
                  legend cell align={left},
                  legend columns = 3,
                  legend style={at={(1.88,1.25)},anchor=north}]
        \addplot+[color=black,solid,mark=o,mark repeat={400}] table[x=t,y=H] {1DNLS_Ham.txt};
        \addplot+[color=blue,mark=none] table[x=t,y=Hs] {1DNLS_Ham.txt};
        \addplot+[color=red,dashed,mark=none,line width=1] table[x=t,y=Hd] {1DNLS_Ham.txt};
        \legend{{High-fidelity},{Static PBDW},{Dyn-PBDW}};
    \nextgroupplot[xlabel={$t$},
                  ylabel={$\eHts(t)$},
                  xlabel style = {yshift=.2cm},
                  ylabel style = {yshift=-.22cm},
                  axis line style = thick,
                  grid=both,
                  grid style = {gray,opacity=0.2},
                  xmin=0, xmax=20,
                  ymin=1e-06, ymax=10^2,
                  ytick={1e-06,1e-04,1e-02,1e+00,1e+02},
                  yticklabels={},
                  extra y ticks={1e-06,1e-04,1e-02,1e+00,1e+02},
                  ymode=log,
                  xlabel style={font=\footnotesize},
                  ylabel style={font=\footnotesize},
                  x tick label style={font=\footnotesize},
                  y tick label style={font=\footnotesize}]
        \addplot+[color=blue,solid,mark=none] table[x=t,y=eHs] {1DNLS_Ham.txt};
        \addplot+[color=red,dashed,mark=none,line width=1] table[x=t,y=eHd] {1DNLS_Ham.txt};
        \nextgroupplot[xlabel={$t$},
                  ylabel={$\dHts(t)$},
                  xlabel style = {yshift=.2cm},
                  ylabel style = {yshift=-.32cm},
                  axis line style = thick,
                  grid=both,
                  grid style = {gray,opacity=0.2},
                  xmin=0, xmax=20,
                  ymin=1e-10, ymax=1e+02,
                  ytick={1e-10,1e-08,1e-06,1e-04,1e-02,1e+00,1e+02},
                  yticklabels={},
                  extra y ticks={1e-10,1e-06,1e-02,1e+02},
                  ymode=log,
                  xlabel style={font=\footnotesize},
                  ylabel style={font=\footnotesize},
                  x tick label style={font=\footnotesize},
                  y tick label style={font=\footnotesize}]
        \addplot+[color=black,solid,mark=o,mark repeat={399}] table[x=t,y=dH] {1DNLS_Ham.txt};
        \addplot+[color=blue,mark=none] table[x=t,y=dHs] {1DNLS_Ham.txt};
        \addplot+[color=red,dashed,mark=none,line width=1] table[x=t,y=dHd] {1DNLS_Ham.txt};
    \end{groupplot}
\end{tikzpicture}
\caption{\footnotesize 1D Schr\"odinger equation. Left: evolution of the Hamiltonian for the test parameter $\theta=(1.0267,0.9867)$. Center: Hamiltonian approximation error \eqref{eq:eHam}. Right: error \eqref{eq:dHam} in the conservation of the Hamiltonian. Comparison between the high-fidelity solution and the reconstructed solution obtained with $m=6$ sensors in the static and dynamic case.}\label{fig:1DNLS_Hamcons}
\end{figure}

Finally, we analyze the conservation of the Hamiltonian associated with the reconstruction. \Cref{fig:1DNLS_Hamcons}, left, shows a comparison between the evolution of the Hamiltonian evaluated at the high-fidelity solution and at the reconstructed solution for the test parameter $\theta=(1.0267,0.9867)$. If the sensors positions are not updated, the Hamiltonian evaluated at the reconstructed solution is close to the target only at the initial stages of the simulation. On the other hand, Dyn-PBDW keeps the error in the preservation of the Hamiltonian below $10^{-2}$ until the final time. We observe that the order of magnitude of both \eqref{eq:eHam} and \eqref{eq:dHam} is comparable to the reconstruction error at all times, in accordance with the analysis of \Cref{sec:Ham_pres}.

In this sections we worked with local measurements of the real and imaginary parts of the complex wave function $\psi$. We point out that, in practice, these quantities might not be available as observables, and what we usually measure instead are quantities related to the modulus $\lvert\psi\rvert=\sqrt{q^2+p^2}$. This is a nonlinear function of $u=(q, p)$, and it is not covered in our current setting, where the assumption on $\ell_i$ being linear plays a crucial role. This assumption can actually be removed in the spirit of \cite{CDMS2022} but we defer this extension to future works.


\subsection{1D shallow water equations}
The one-dimensional, nonlinear shallow water equations read
\begin{equation*}
    \begin{cases}
        h_t+(h\Phi_x)_x=0 & \quad x\in\Sdom, \; t\in[0,T], \\
        \Phi_t+\frac{1}{2}(\Phi_x)^2+h=0 & \quad x\in\Sdom, \; t\in[0,T],
    \end{cases}
\end{equation*}
where $h$ is the height of the free surface and $\Phi$ is the velocity potential. We complement the problem with periodic boundary conditions. The initial condition is $h(x,0)=1+\alpha\exp(-\beta x^2)$  where $\theta=(\alpha,\beta)\in\Theta=[1/10,1/7]\times[2/10,15/10]$, while $\Phi(x,0)=0$. The Hamiltonian for this problem is
\begin{equation*}
    \mathcal{H}(u)=\frac{1}{2}\int_\Sdom h(\Phi_x^2+h)\,dx, \qquad u=(h,\Phi).
\end{equation*}
Notice that, since the parameters $\alpha$ and $\beta$ only appear in the initial condition, the Hamiltonian $\Hcal$ is independent of $\theta$.

\begin{figure}[ht]
    \centering
    \includegraphics[width=0.95\textwidth]{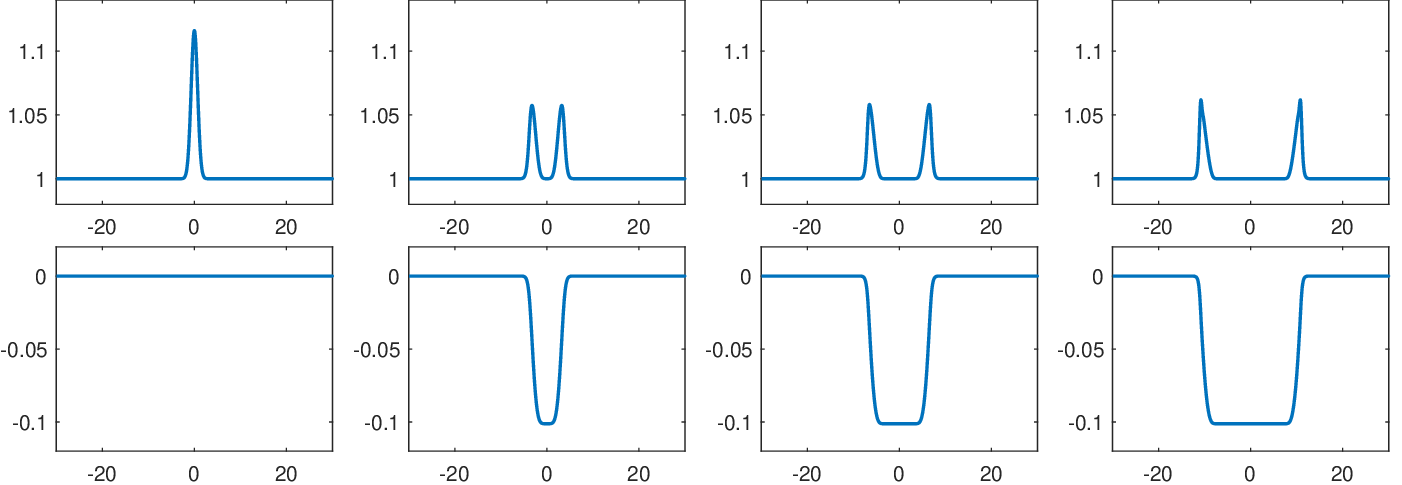}
    \caption{\footnotesize 1D shallow water equations. Depth $h(x,t)$ (above) and velocity potential $\Phi(x,t)$ (below) at times $t=0$, $t=3$, $t=6$ and $t=10$ for $\theta=(0.1161,1.0125)$.}\label{fig:1DSWE_fullsol}
\end{figure}

In the numerical experiments we set $L_x=30$, $N_x=1000$, $T=10$ and $N_t=10000$. The evolution of $h$ and $\Phi$ is shown in \Cref{fig:1DSWE_fullsol} for the test parameter $\theta=(0.1161,1.0125)$: the initial profile splits into two crests that move in opposite directions.

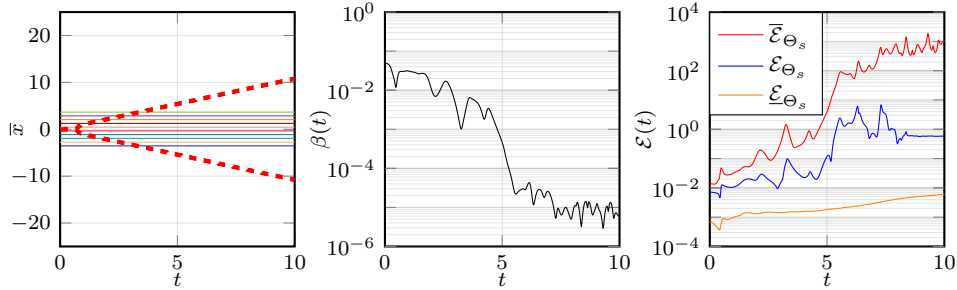
\begin{figure}[ht]
\centering
\begin{tikzpicture}
    \begin{groupplot}[
      group style={group size=3 by 1, horizontal sep=1.2cm},
      width=4.7cm, height=4.7cm
    ]
    \nextgroupplot[xlabel={$t$},
                  ylabel={$\xs$},
                  xlabel style = {yshift=.2cm},
                  ylabel style = {yshift=-.4cm},
                  xtick={0,5,10},
                  ytick={-20,-10,0,10,20},
                  axis line style = thick,
                  grid=both,
                  grid style = {gray,opacity=0.2},
                  xmin=0, xmax=10,
                  ymin=-25, ymax=25,
                  xlabel style={font=\footnotesize},
                  ylabel style={font=\footnotesize},
                  x tick label style={font=\footnotesize},
                  y tick label style={font=\footnotesize}]
        \addplot+[color=Blue,mark=none] table[x=t,y=s1x] {1DSWE_noadapt_pointwise.txt};
        \addplot+[color=Dandelion,mark=none] table[x=t,y=s2x] {1DSWE_noadapt_pointwise.txt};
        \addplot+[color=Cerulean,mark=none] table[x=t,y=s3x] {1DSWE_noadapt_pointwise.txt};
        \addplot+[color=PineGreen,mark=none] table[x=t,y=s4x] {1DSWE_noadapt_pointwise.txt};
        \addplot+[color=OrangeRed,mark=none] table[x=t,y=s5x] {1DSWE_noadapt_pointwise.txt};
        \addplot+[color=Tan,solid,mark=none] table[x=t,y=s6x] {1DSWE_noadapt_pointwise.txt};
        \addplot+[color=Sepia,solid,mark=none] table[x=t,y=s7x] {1DSWE_noadapt_pointwise.txt};
        \addplot+[color=Orange,solid,mark=none] table[x=t,y=s8x] {1DSWE_noadapt_pointwise.txt};
        \addplot+[color=Fuchsia,solid,mark=none] table[x=t,y=s9x] {1DSWE_noadapt_pointwise.txt};
        \addplot+[color=LimeGreen,solid,mark=none] table[x=t,y=s10x] {1DSWE_noadapt_pointwise.txt};
        \addplot+[color=red,dashed,mark=none,line width=1.5pt] table[x=t,y=bump1x] {1DSWE_noadapt_pointwise.txt};
        \addplot+[color=red,dashed,mark=none,line width=1.5pt] table[x=t,y=bump2x] {1DSWE_noadapt_pointwise.txt};
    \nextgroupplot[xlabel={$t$},
                  ylabel={$\beta(t)$},
                  xlabel style = {yshift=.2cm},
                  ylabel style = {yshift=-.3cm},
                  xtick={0,5,10},
                  ytick={1e-06,1e-05,1e-04,1e-03,1e-02,1e-01,1},
                  yticklabels={},
                  extra y ticks={1e-06,1e-04,1e-02,1e+00},
                  yminorticks = false,
                  axis line style = thick,
                  grid=both,
                  ymode=log,
                  grid style = {gray,opacity=0.2},
                  xmin=0, xmax=10,
                  ymin=1e-06, ymax=1,
                  xlabel style={font=\footnotesize},
                  ylabel style={font=\footnotesize},
                  x tick label style={font=\footnotesize},
                  y tick label style={font=\footnotesize}]
        \addplot+[color=black,mark=none] table[x=t,y=beta] {1DSWE_noadapt_pointwise.txt};
    \nextgroupplot[xlabel={$t$},
                  ylabel={$\Ecal(t)$},
                  xlabel style = {yshift=.2cm},
                  ylabel style = {yshift=-.3cm},
                  axis line style = thick,
                  grid=both,
                  ytick={1e-04,1e-03,1e-02,1e-01,1e+00,1e+01,1e+02,1e+03,1e+04},
                  yticklabels={},
                  extra y ticks={1e-04,1e-02,1e+00,1e+02,1e+04},
                  yminorticks = false,
                  xtick={0,5,10},
                  ymode=log,
                  grid style = {gray,opacity=0.2},
                  xmin=0, xmax=10,
                  ymin=1e-04, ymax=1e+04,
                  xlabel style={font=\footnotesize},
                  ylabel style={font=\footnotesize},
                  x tick label style={font=\footnotesize},
                  y tick label style={font=\footnotesize},
                  legend style={font=\footnotesize},
                  legend cell align={left},
                  legend columns = 1,
                  legend style={at={(0.24,1)},anchor=north}]
        \addplot+[color=red,mark=none] table[x=t,y=eb] {1DSWE_noadapt_pointwise.txt};
        \addplot+[color=blue,mark=none] table[x=t,y=ese] {1DSWE_noadapt_pointwise.txt};
        \addplot+[color=orange,mark=none] table[x=t,y=ep] {1DSWE_noadapt_pointwise.txt};
        \legend{{$\Eboundts$},{$\Esets$},{$\Eprojts$}};
    \end{groupplot}
\end{tikzpicture}
\caption{\footnotesize 1D shallow water equations. Static local averages of width $\sigma=0.1$, $m=10$. Left: $x$ coordinate of the sensors over time; the red dashed lines denote the position of the main humps of the high-fidelity solution $h(t,\theta)$ for $\theta=(0.1161,1.0125)$. Center: evolution of the inf-sup constant $\beta$. Right: evolution of the numerical errors with $\lvert\Thtest\rvert=16$ test parameters.}\label{fig:1DSWE_noadapt_pointwise}
\end{figure}

Concerning the evolution of the approximation space, we set $2n=12$ and $\lvert\Thtrain\rvert=5^2=25$. State estimation is performed with $m=10$ sensors for $\lvert\Thtest\rvert=4^2=16$ test parameters. At the initial time the sensors are equispaced in the sub-interval $[-3.5,3.5]$, around the main hump. Results obtained with static measurements are depicted in \Cref{fig:1DSWE_noadapt_pointwise}. As time evolves, $\beta$ decays until the final value of approximately $10^{-6}$. Moreover, the state estimation error $\Esets$ is at least one order of magnitude larger than the best approximation error $\Eprojts$, even at the initial time. This implies that the initial placement of the sensors around the main profile, although reasonable, does not correspond to the best possible configuration.

This can be addressed by Dyn-PBDW as shown in \Cref{fig:1DSWE_adapt_pointwise}. We observe that running \Cref{alg:grad_desc} at the initial time improves the stability threshold obtained with the initial arbitrary guess of the sensors positions, as seen by comparing the values of $\beta(0)$ in the static and dynamic case. Moreover, the sensors are able to follow both profiles that characterize the high-fidelity solution, and the value of $\beta$ remains above $10^{-1}$ in the whole temporal window.

\begin{figure}[ht]
\centering
\begin{tikzpicture}
    \begin{groupplot}[
      group style={group size=3 by 1, horizontal sep=1.2cm},
      width=4.7cm, height=4.7cm
    ]
    \nextgroupplot[xlabel={$t$},
                  ylabel={$\xs$},
                  xlabel style = {yshift=.2cm},
                  ylabel style = {yshift=-.4cm},
                  xtick={0,5,10},
                  ytick={-20,-10,0,10,20},
                  axis line style = thick,
                  grid=both,
                  grid style = {gray,opacity=0.2},
                  xmin=0, xmax=10,
                  ymin=-25, ymax=25,
                  xlabel style={font=\footnotesize},
                  ylabel style={font=\footnotesize},
                  x tick label style={font=\footnotesize},
                  y tick label style={font=\footnotesize}]
        \addplot+[color=Blue,mark=none] table[x=t,y=s1x] {1DSWE_adapt_pointwise.txt};
        \addplot+[color=Dandelion,mark=none] table[x=t,y=s2x] {1DSWE_adapt_pointwise.txt};
        \addplot+[color=Cerulean,mark=none] table[x=t,y=s3x] {1DSWE_adapt_pointwise.txt};
        \addplot+[color=PineGreen,mark=none] table[x=t,y=s4x] {1DSWE_adapt_pointwise.txt};
        \addplot+[color=OrangeRed,mark=none] table[x=t,y=s5x] {1DSWE_adapt_pointwise.txt};
        \addplot+[color=Tan,solid,mark=none] table[x=t,y=s6x] {1DSWE_adapt_pointwise.txt};
        \addplot+[color=Sepia,solid,mark=none] table[x=t,y=s7x] {1DSWE_adapt_pointwise.txt};
        \addplot+[color=Orange,solid,mark=none] table[x=t,y=s8x] {1DSWE_adapt_pointwise.txt};
        \addplot+[color=Fuchsia,solid,mark=none] table[x=t,y=s9x] {1DSWE_adapt_pointwise.txt};
        \addplot+[color=LimeGreen,solid,mark=none] table[x=t,y=s10x] {1DSWE_adapt_pointwise.txt};
        \addplot+[color=red,dashed,mark=none,line width=1.5pt] table[x=t,y=bump1x] {1DSWE_adapt_pointwise.txt};
        \addplot+[color=red,dashed,mark=none,line width=1.5pt] table[x=t,y=bump2x] {1DSWE_adapt_pointwise.txt};
    \nextgroupplot[xlabel={$t$},
                  ylabel={$\beta(t)$},
                  xlabel style = {yshift=.2cm},
                  ylabel style = {yshift=-.3cm},
                  axis line style = thick,
                  grid=both,
                  xtick={0,5,10},
                  ymode=log,
                  grid style = {gray,opacity=0.2},
                  xmin=0, xmax=10,
                  ymin=1e-06, ymax=1e+00,
                  ytick={1e-06,1e-05,1e-04,1e-03,1e-02,1e-01,1e+00},
                  yticklabels = {},
                  extra y ticks = {1e-06,1e-04,1e-02,1e+00},
                  yminorticks = false,
                  xlabel style={font=\footnotesize},
                  ylabel style={font=\footnotesize},
                  x tick label style={font=\footnotesize},
                  y tick label style={font=\footnotesize}]
        \addplot+[color=black,mark=none] table[x=t,y=beta] {1DSWE_adapt_pointwise.txt};
    \nextgroupplot[xlabel={$t$},
                  ylabel={$\Ecal(t)$},
                  xlabel style = {yshift=.2cm},
                  ylabel style = {yshift=-.3cm},
                  axis line style = thick,
                  ytick={1e-03,1e+00,1e+03},
                  xtick={0,5,10},
                  grid=both,
                  ymode=log,
                  grid style = {gray,opacity=0.2},
                  xmin=0, xmax=10,
                  ymin=1e-04, ymax=1e+04,
                  ytick={1e-04,1e-03,1e-02,1e-01,1e+00,1e+01,1e+02,1e+03,1e+04},
                  yticklabels = {},
                  extra y ticks = {1e-04,1e-02,1e+00,1e+02,1e+04},
                  yminorticks = false,
                  xlabel style={font=\footnotesize},
                  ylabel style={font=\footnotesize},
                  x tick label style={font=\footnotesize},
                  y tick label style={font=\footnotesize},
                  legend style={font=\footnotesize},
                  legend cell align={left},
                  legend columns = 1,
                  legend style={at={(0.27,0.98)},anchor=north}]
        \addplot+[color=red,mark=none] table[x=t,y=eb] {1DSWE_adapt_pointwise.txt};
        \addplot+[color=blue,mark=none] table[x=t,y=ese] {1DSWE_adapt_pointwise.txt};
        \addplot+[color=orange,mark=none] table[x=t,y=ep] {1DSWE_adapt_pointwise.txt};
        \legend{{$\Eboundts$},{$\Esets$},{$\Eprojts$}};
    \end{groupplot}
\end{tikzpicture}
\caption{\footnotesize 1D shallow water equations. Dynamic local averages of width $\sigma=0.1$, $m=10$. Left: $x$ coordinate of the sensors over time; the red dashed lines denote the position of the main bumps of the high-fidelity solution $h(t,\theta)$ for $\theta=(0.1161,1.0125)$. Center: evolution of the inf-sup constant $\beta$. Right: evolution of the numerical errors with $\lvert\Thtest\rvert=16$ test parameters.}\label{fig:1DSWE_adapt_pointwise}
\end{figure}

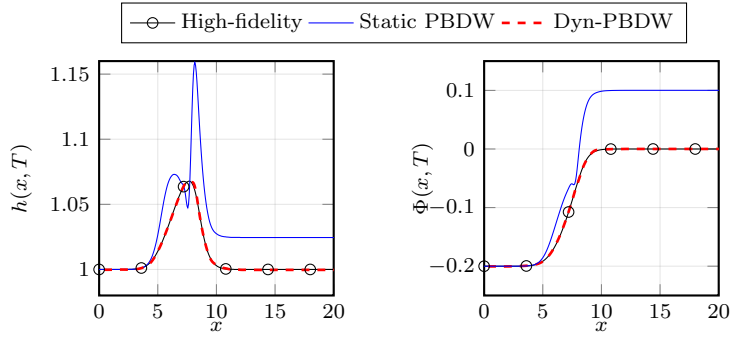
\begin{figure}[ht]
\centering
\begin{tikzpicture}
    \begin{groupplot}[
      group style={group size=2 by 1, horizontal sep=2cm},
      width=4.7cm, height=4.7cm
    ]
    \nextgroupplot[xlabel={$x$},
                  ylabel={$h(x,T)$},
                  xlabel style = {yshift=.2cm},
                  axis line style = thick,
                  grid=both,
                  grid style = {gray,opacity=0.2},
                  xmin=0, xmax=20,
                  ymin=0.98, ymax=1.16,
                  xlabel style={font=\footnotesize},
                  ylabel style={font=\footnotesize},
                  x tick label style={font=\footnotesize},
                  y tick label style={font=\footnotesize},
                  legend style={font=\footnotesize},
                  legend cell align={left},
                  legend columns = 3,
                  legend style={at={(1.3,1.25)},anchor=north}]
        \addplot+[color=black,mark=o,mark repeat={60}] table[x=xgrid,y=hex] {1DSWE_recsol_a_vs_na.txt};
        \addplot+[color=blue,mark=none] table[x=xgrid,y=hna] {1DSWE_recsol_a_vs_na.txt};
        \addplot+[color=red,dashed,mark=none,line width=1] table[x=xgrid,y=ha] {1DSWE_recsol_a_vs_na.txt};
        \legend{{High-fidelity},{Static PBDW},{Dyn-PBDW}};
    \nextgroupplot[xlabel={$x$},
                  ylabel={$\Phi(x,T)$},
                  xlabel style = {yshift=.2cm},
                  ylabel style = {yshift=-.3cm},
                  axis line style = thick,
                  grid=both,
                  grid style = {gray,opacity=0.2},
                  xmin=0, xmax=20,
                  ymin=-0.25, ymax=0.15,
                  xlabel style={font=\footnotesize},
                  ylabel style={font=\footnotesize},
                  x tick label style={font=\footnotesize},
                  y tick label style={font=\footnotesize},
                  ]
        \addplot+[color=black,mark=o,mark repeat={60}] table[x=xgrid,y=phiex] {1DSWE_recsol_a_vs_na.txt};
        \addplot+[color=red,dashed,mark=none,line width=1] table[x=xgrid,y=phia] {1DSWE_recsol_a_vs_na.txt};
        \addplot+[color=blue,mark=none] table[x=xgrid,y=phina] {1DSWE_recsol_a_vs_na.txt};
    \end{groupplot}
\end{tikzpicture}
\caption{\footnotesize 1D shallow water equations. Numerical solution of the high-fidelity model corresponding to $\theta=(0.1375,0.3625)$ and comparison with the statically and dynamically reconstructed solutions at time $t=7$. For the sake of visualization, only the interval $[0,20]\subset\Sdom$ is shown.}\label{fig:1DSWE_recsol}
\end{figure}

We show the reconstructed solutions at time $t=7$ in \Cref{fig:1DSWE_recsol} for one test parameter, and a comparison with the high-fidelity solution. For the sake of visualization, only one of the two crests is shown. If the sensors are not moved, the reconstruction fails to correctly identify its amplitude and height, as well as its velocity. On the other hand, the reconstruction obtained with Dyn-PBDW is close to the ground truth for both $h$ and $\Phi$.

We conclude this section by analyzing the evolution of the Hamiltonian evaluated at the reconstructed solution. As shown in \Cref{fig:1DSWE_Hamcons}, left, it is not possible to accurately approximate the high-fidelity Hamiltonian without adapting the positions of the sensors. In fact, \Cref{fig:1DSWE_Hamcons}, center, shows that the absolute approximation error in this case can grow up to over $10^1$ in the worst case scenario, while Dyn-PBDW is able to keep this quantity below $10^{-1}$. Finally, we report in \Cref{fig:1DSWE_Hamcons}, right, the preservation of the Hamiltonian with respect to its value at the initial time. We observe here that the errors in the dynamic case exhibit a moderate growth. This is to be expected, since the reconstruction step is not guaranteed to preserve the Hamiltonian. Moreover, as observed in \Cref{sec:Ham_pres}, these errors are related to the reconstruction error which, in turn, is bounded from below by the projection error onto the approximation space. As seen before, the latter is also around $10^{-2}$ at the final time, which explains the behaviour shown in \Cref{fig:1DSWE_Hamcons}. We expect the quality of the reconstruction to improve if the dynamical approximation discussed in \Cref{sec:Ham_dynamics} is combined with a rank-adaptive strategy to allow the dimension of the approximation space to increase over time and prevent the projection error from growing. This is left for future work.

\begin{figure}[ht]
\centering
\begin{tikzpicture}
    \begin{groupplot}[
      group style={group size=3 by 1, horizontal sep=1.2cm},
      width=4.7cm, height=4.7cm
    ]
    \nextgroupplot[xlabel={$t$},
                  ylabel={$\Hcal(\mathfrak{u}(t))$},
                  xlabel style = {yshift=.2cm},
                  ylabel style = {yshift=-.2cm},
                  axis line style = thick,
                  grid=both,
                  grid style = {gray,opacity=0.2},
                  xmin=0, xmax=10,
                  ymin=20, ymax=65,
                  xlabel style={font=\footnotesize},
                  ylabel style={font=\footnotesize},
                  x tick label style={font=\footnotesize},
                  y tick label style={font=\footnotesize},
                  legend style={font=\footnotesize},
                  legend cell align={left},
                  legend columns = 3,
                  legend style={at={(1.88,1.25)},anchor=north}]
        \addplot+[color=black,mark=o,mark repeat={199}] table[x=t,y=H] {1DSWE_Ham.txt};
        \addplot+[color=blue,mark=none] table[x=t,y=Hs] {1DSWE_Ham.txt};
        \addplot+[color=red,dashed,mark=none,line width=1] table[x=t,y=Hd] {1DSWE_Ham.txt};
        \legend{{High-fidelity},{Static PBDW},{Dyn-PBDW}};
    \nextgroupplot[xlabel={$t$},
                  ylabel={$\eHts(t)$},
                  xlabel style = {yshift=.2cm},
                  ylabel style = {yshift=-.22cm},
                  axis line style = thick,
                  grid=both,
                  grid style = {gray,opacity=0.2},
                  xmin=0, xmax=10,
                  ymin=1e-04, ymax=1e+02,
                  ytick={1e-04,1e-03,1e-02,1e-01,1e+00,1e+01,1e+02},
                  yticklabels={},
                  extra y ticks={1e-04,1e-02,1e+00,1e+02},
                  yminorticks=false,
                  ymode=log,
                  xlabel style={font=\footnotesize},
                  ylabel style={font=\footnotesize},
                  x tick label style={font=\footnotesize},
                  y tick label style={font=\footnotesize}]
        \addplot+[color=blue,solid,mark=none] table[x=t,y=eHs] {1DSWE_Ham.txt};
        \addplot+[color=red,dashed,mark=none,line width=1] table[x=t,y=eHd] {1DSWE_Ham.txt};
    \nextgroupplot[xlabel={$t$},
                  ylabel={$\dHts(t)$},
                  xlabel style = {yshift=.2cm},
                  ylabel style = {yshift=-.32cm},
                  axis line style = thick,
                  grid=both,
                  grid style = {gray,opacity=0.2},
                  xmin=0, xmax=10,
                  ymin=1e-13, ymax=1e+03,
                  ytick={1e-13,1e-11,1e-09,1e-07,1e-05,1e-03,1e-01,1e+01,1e+03},
                  yticklabels={},
                  extra y ticks={1e-13,1e-09,1e-05,1e-01,1e+03},
                  yminorticks=false,
                  ymode=log,
                  xlabel style={font=\footnotesize},
                  ylabel style={font=\footnotesize},
                  x tick label style={font=\footnotesize},
                  y tick label style={font=\footnotesize}]
        \addplot+[color=black,solid,mark=o,mark repeat={199}] table[x=t,y=dH] {1DSWE_Ham.txt};
        \addplot+[color=blue,mark=none] table[x=t,y=dHs] {1DSWE_Ham.txt};
        \addplot+[color=red,dashed,mark=none,line width=1] table[x=t,y=dHd] {1DSWE_Ham.txt};
    \end{groupplot}
\end{tikzpicture}
\caption{\footnotesize 1D shallow water equations. Left: evolution of the Hamiltonian for the test parameter $\theta=(0.1161,0.3625)$. Center: Hamiltonian approximation error \eqref{eq:eHam}. Right: error \eqref{eq:dHam} in the conservation of the Hamiltonian. Comparison between the high-fidelity solution and the reconstructed solution obtained with $m=10$ sensors in the static and dynamic case.}\label{fig:1DSWE_Hamcons}
\end{figure}


\subsection{2D shallow water equations}
As a final numerical test, we consider the two-dimensional nonlinear shallow water equations
\begin{equation*}
    \begin{cases}
        h_t+\nabla\cdot(h\nabla\Phi)=0 & x\in\Sdom \quad t\in[0,T] \\
        \Phi_t+\frac{1}{2}\lvert\nabla\Phi\rvert^2+h=0 & x\in\Sdom \quad t\in[0,T]
    \end{cases}
\end{equation*}
with periodic boundary conditions. The initial condition is $\Phi(x,0)=0$ and $h(x,0)=1+\alpha\exp{(-\beta\lVert x\rVert^2})$, where $\theta=(\alpha,\beta)\in\Theta=[1/5,1/2]\times[11/10,17/10]$. We set $L_x=L_y=8$, and $N_x=N_y=50$ intervals in each direction. The final time is $T=10$ and $N_t=10000$ uniform steps are used for time integration. 

\begin{figure}[h]
    \centering
    \includegraphics[width=0.95\textwidth]{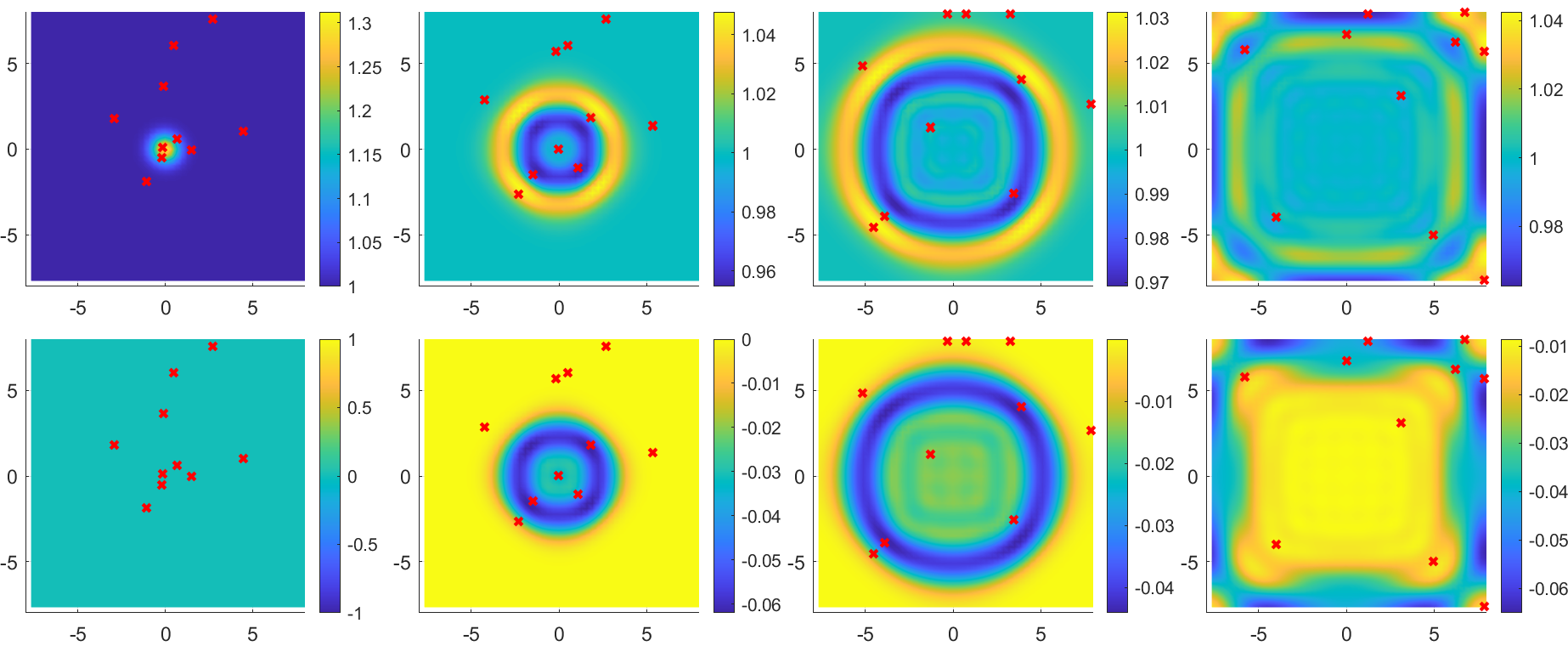}
    \caption{\footnotesize 2D shallow water equations. Depth $h(x,y,t)$ (above) and velocity potential $\Phi(x,y,t)$ (below) for the test parameter $\theta=(0.3125,1.6250)$ at times $t=0$, $t=3$, $t=6$ and $t=10$. The red crosses denote the locations of the sensors.}
    \label{fig:2DSWE_senspos}
\end{figure}

The approximation space, of dimension $2n=12$, is obtained with $\lvert\Thtrain\rvert=5^2=25$ uniformly selected parameters in $\Theta$. We then consider $m=10$ sensors for state estimation. The initial positions of the sensors are randomly selected in the square $[-0.8,0.8]^2$, around the initial profile. \Cref{fig:2DSWE_senspos} shows the location of the measurements computed with Dyn-PBDW.

\begin{figure}[ht]
\centering
\begin{tikzpicture}
    \begin{groupplot}[
      group style={group size=2 by 1, horizontal sep=2cm},
      width=4.7cm, height=4.7cm]
    \nextgroupplot[xlabel={$t$},
                  ylabel={$\beta(t)$},
                  xlabel style = {yshift=.2cm},
                  ylabel style = {yshift=-.2cm},
                  axis line style = thick,
                  grid=both,
                  ymode=log,
                  grid style = {gray,opacity=0.2},
                  xmin=0, xmax=10,
                  ymin=1e-08, ymax=1e+00,
                  ytick={1e-08,1e-07,1e-06,1e-05,1e-04,1e-03,1e-02,1e-01,1e+00},
                  yticklabels={},
                  extra y ticks={1e-08,1e-06,1e-04,1e-02,1e+00},
                  yminorticks=false,
                  xlabel style={font=\footnotesize},
                  ylabel style={font=\footnotesize},
                  x tick label style={font=\footnotesize},
                  y tick label style={font=\footnotesize},
                  legend style={font=\footnotesize},
                  legend cell align={left},
                  legend columns = 2,
                  legend style={at={(0.5,1.22)},anchor=north}]
        \addplot+[color=black,solid,mark=none] table[x=t,y=betana]{2DSWE_sigma0p1_errbeta.txt};
        \addplot+[color=black,dashed,mark=none] table[x=t,y=betaa] {2DSWE_sigma0p1_errbeta.txt};
        \legend{{Static PBDW},{Dyn-PBDW}};
    \nextgroupplot[xlabel={$t$},
                  ylabel={$\Ecal(t)$},
                  xlabel style = {yshift=.2cm},
                  ylabel style = {yshift=-.3cm},
                  axis line style = thick,
                  grid=both,
                  ymode=log,
                  grid style = {gray,opacity=0.2},
                  xmin=0, xmax=10,
                  ymin=1e-06, ymax=1e+04,
                  ytick={1e-06,1e-05,1e-04,1e-03,1e-02,1e-01,1e+00,1e+01,1e+02,1e+03,1e+04},
                  yticklabels={},
                  extra y ticks={1e-06,1e-03,1e+00,1e+03},
                  yminorticks = false,
                  xlabel style={font=\footnotesize},
                  ylabel style={font=\footnotesize},
                  x tick label style={font=\footnotesize},
                  y tick label style={font=\footnotesize},
                  legend style={font=\footnotesize},
                  legend cell align={left},
                  legend columns = 1,
                  legend style={at={(1.3,0.85)},anchor=north}]
        \addplot+[color=red,solid,mark=none] table[x=t,y=ebna] {2DSWE_sigma0p1_errbeta.txt};
        \addplot+[color=blue,solid,mark=none] table[x=t,y=esena] {2DSWE_sigma0p1_errbeta.txt};
        \addplot+[color=red,dashed,mark=none] table[x=t,y=eba] {2DSWE_sigma0p1_errbeta.txt};
        \addplot+[color=blue,dashed,mark=none] table[x=t,y=esea] {2DSWE_sigma0p1_errbeta.txt};
        \addplot+[color=orange,solid,mark=none] table[x=t,y=ep] {2DSWE_sigma0p1_errbeta.txt};
        \legend{{$\Eboundts^{st}$},{$\Esets^{st}$},{$\Eboundts^{dyn}$},{$\Esets^{dyn}$},{$\Eprojts$}};
    \end{groupplot}
\end{tikzpicture}
\caption{\footnotesize 2D shallow water equations. Local averages of width $\sigma=0.1$, $m=10$. Left: evolution of the inf-sup constant $\beta$ over time. Right: evolution of the numerical errors with $\lvert\Thtest\rvert=16$ test parameters.}\label{fig:2DSWE_sigma0p1_errbeta}
\end{figure}
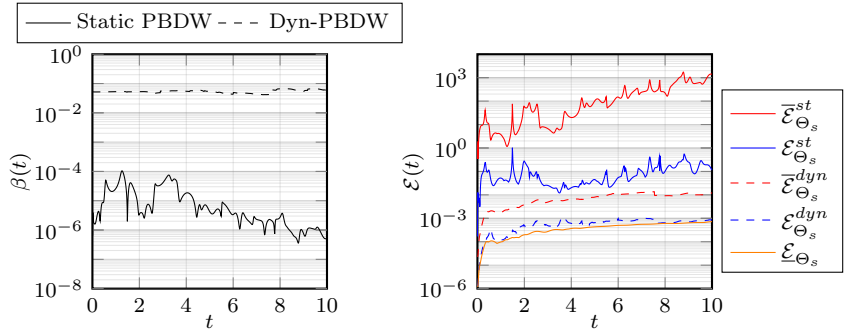

The evolution of $\beta$ and of the state estimation errors obtained with $\lvert\Thtest\rvert=16$ test parameters are reported in \Cref{fig:2DSWE_sigma0p1_errbeta}. As expected, we notice a moderate decay of the value of $\beta$ when the sensor locations are not updated, while Dyn-PBDW allows to keep it approximately constant (\Cref{fig:2DSWE_sigma0p1_errbeta}, left). Furthermore, as in the one-dimensional case, it is possible to improve the stability threshold obtained with the initial random guess, as we notice by comparing the values of $\beta(0)$ in the static and dynamic case. This can be useful in situations where it is not clear how to optimally place the sensors at the initial time or when the initial condition of the problem is not known at all.

\begin{figure}[ht]
\centering
\begin{tikzpicture}
    \begin{groupplot}[
      group style={group size=3 by 1, horizontal sep=1.2cm},
      width=4.7cm, height=4.7cm
    ]
    \nextgroupplot[xlabel={$t$},
                  ylabel={$\Hcal(\mathfrak{u}(t))$},
                  xlabel style = {yshift=.2cm},
                  ylabel style = {yshift=-.5cm},
                  axis line style = thick,
                  grid=both,
                  grid style = {gray,opacity=0.2},
                  xmin=0, xmax=10,
                  ymin=385, ymax=415,
                  ytick={385,400,415},
                  yticklabels={},
                  extra y ticks={385,415},
                  xlabel style={font=\footnotesize},
                  ylabel style={font=\footnotesize},
                  x tick label style={font=\footnotesize},
                  y tick label style={font=\footnotesize},
                  legend style={font=\footnotesize},
                  legend cell align={left},
                  legend columns = 3,
                  legend style={at={(1.88,1.25)},anchor=north}]
        \addplot+[color=black,mark=o,mark repeat={199}] table[x=t,y=H] {2DSWE_Ham.txt};
        \addplot+[color=blue,mark=none] table[x=t,y=Hs] {2DSWE_Ham.txt};
        \addplot+[color=red,dashed,mark=none,line width=1.5] table[x=t,y=Hd] {2DSWE_Ham.txt};
        \legend{{High-fidelity},{Static PBDW},{Dyn-PBDW}};
    \nextgroupplot[xlabel={$t$},
                  ylabel={$\eHts(t)$},
                  xlabel style = {yshift=.2cm},
                  ylabel style = {yshift=-.22cm},
                  axis line style = thick,
                  grid=both,
                  grid style = {gray,opacity=0.2},
                  xmin=0, xmax=10,
                  ymin=1e-06, ymax=1e+04,
                  ytick={1e-06,1e-05,1e-04,1e-03,1e-02,1e-01,1e+00,1e+01,1e+02,1e+03,1e+04},
                  yticklabels={},
                  extra y ticks={1e-06,1e-04,1e-02,1e+00,1e+02,1e+04},
                  yminorticks = false,
                  ymode=log,
                  xlabel style={font=\footnotesize},
                  ylabel style={font=\footnotesize},
                  x tick label style={font=\footnotesize},
                  y tick label style={font=\footnotesize}]
        \addplot+[color=blue,solid,mark=none] table[x=t,y=eHs] {2DSWE_Ham.txt};
        \addplot+[color=red,dashed,mark=none,line width=1.5] table[x=t,y=eHd] {2DSWE_Ham.txt};
    \nextgroupplot[xlabel={$t$},
                  ylabel={$\dHts(t)$},
                  xlabel style = {yshift=.2cm},
                  ylabel style = {yshift=-.32cm},
                  axis line style = thick,
                  grid=both,
                  grid style = {gray,opacity=0.2},
                  xmin=0, xmax=10,
                  ymin=1e-12, ymax=1e+04,
                  ytick={1e-12,1e-10,1e-08,1e-06,1e-04,1e-02,1e+00,1e+02,1e+04},
                  yticklabels={},
                  extra y ticks={1e-12,1e-08,1e-04,1e+00,1e+04},
                  ymode=log,
                  xlabel style={font=\footnotesize},
                  ylabel style={font=\footnotesize},
                  x tick label style={font=\footnotesize},
                  y tick label style={font=\footnotesize}]
        \addplot+[color=black,solid,mark=o,mark repeat={199}] table[x=t,y=dH] {2DSWE_Ham.txt};
        \addplot+[color=blue,mark=none] table[x=t,y=dHs] {2DSWE_Ham.txt};
        \addplot+[color=red,dashed,mark=none,line width=1.5] table[x=t,y=dHd] {2DSWE_Ham.txt};
    \end{groupplot}
\end{tikzpicture}
\caption{\footnotesize 2D shallow water equations. Left: evolution of the Hamiltonian for the test parameter $\theta=(0.3125,1.6250)$. Center: Hamiltonian approximation error \eqref{eq:eHam}. Right: error \eqref{eq:dHam} in the conservation of the Hamiltonian. Comparison between the high-fidelity solution and the reconstructed solution obtained with $m=10$ sensors in the static and dynamic case.}\label{fig:2DSWE_Hamcons}
\end{figure}
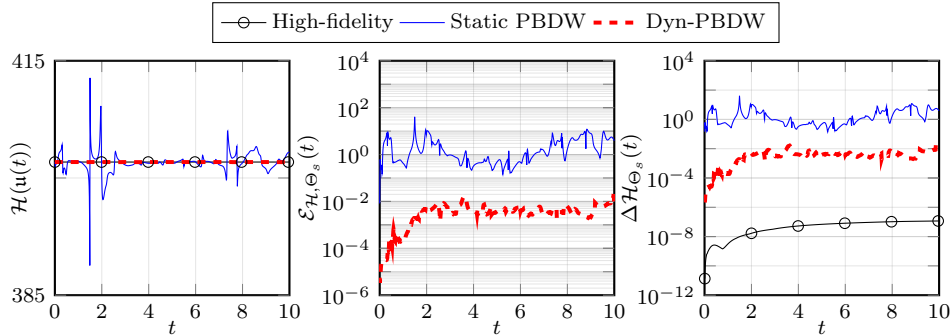

Finally, we show in \Cref{fig:2DSWE_Hamcons} the evolution of the Hamiltonian for one test parameter and the Hamiltonian approximation errors. Similar conclusions as in the one dimensional case hold. In particular, the errors in the dynamic case are always at least two orders of magnitude lower than those achieved without updating the positions of the sensors.

\section{Conclusions}
We have developed the theoretical foundations of Dyn-PBDW, a filtering algorithm leveraging dynamical approximation. We have introduced a dynamical sensor placement strategy which plays a crucial role to guarantee a good reconstruction quality in transport-dominated problems. Beyond the concrete details of the Dyn-PBDW algorithm, the work presents three important conceptual contributions which we hope will inspire future works:
\begin{itemize}
    \item We have given a rigorous theoretical framework to analyze dynamical sensor placement strategies;
    \item We have provided rigorous arguments in favor of working with dynamical reduced models for data assimilation and filtering;
    \item As a simple by-product of our construction, it is possible to incorporate relevant physical properties in the reconstructed solutions if the dynamical approximate models are appropriately built. This is illustrated by the property of symplecticity in our example with Hamiltonian problems.
\end{itemize}
Among the several possible extensions of the work, we may mention:
\begin{itemize}
\item Extending the theory to work with nonlinear observations;
\item Adding information on the reconstructed flow into the construction of the approximation;
\item Dynamically adapting the probability distribution over the parameter space $\Theta$ in the spirit of Bayesian approaches.
\end{itemize}

\appendix

\section{Practical computation of \texorpdfstring{$\nabla_\xs\beta^2$}{gradb} for Hamiltonian systems}\label{app:gradbeta_Ham}
In this appendix, we specialize to the case of Hamiltonian systems the expression of $\nabla_\xs\beta^2$ derived in \Cref{sec:dyn_sens_plac}. As in \Cref{sec:Ham_dynamics}, let $\Sdom\subset\bR^d$ with $d\geq 1$, and let $\Vqp$ be a real Hilbert space on $\Sdom$, with inner product $\inprodV{\cdot}{\cdot}_{\Vqp}$. The dual space $V^\prime$ of $V=\Vqp\times\Vqp$ is isomorphic to $\Vqp^\prime\times\Vqp^\prime$. Thus, we may identify every $\ell\in V^\prime$ with the pair $(\ell^q,\ell^p)$, where $\ell^q,\ell^p\in\Vqp^\prime$ and $\ell(u):=\ell^q(u^q)+\ell^p(u^p)$ for all $u=(u^q,u^p)\in V$. If $\omega^q$ and $\omega^p$ are the Riesz representers of $\ell^q$ and $\ell^p$ in $\Vqp$, respectively, then $\omega:=(\omega^q,\omega^p)$ is the Riesz representer of $\ell$ in $V$, since
\begin{equation*}
    \inprodV{\omega}{u}=\inprodV{\omega^q}{u^q}_{\Vqp}+\inprodV{\omega^p}{u^p}_{\Vqp}=\ell^q(u^q)+\ell^p(u^p)=\ell(u).
\end{equation*}

For $u\in V$, we consider $2m$ linear measurements $z_j:=\ell_j(u)$, $j=1,\dots,2m$, where $\ell_j^p=0$ for $j=1,\dots,m$ and $\ell_j^q=0$ for $j=m+1,\dots,2m$. In other words,
\begin{equation}\label{eq:liHam}
    z_j=\ell_j(u)=
    \begin{cases}
        \ell_j^q(u^q) & j=1,\dots,m \\
        \ell_j^p(u^p) & j=m+1,\dots,2m
    \end{cases}.
\end{equation}
The observation space is then $\Wtwom:=\text{span}\{\omega_1,\dots,\omega_{2m}\}$, where the Riesz representer $\omega_j\in V$ of $\ell_j$ is
\begin{equation*}
    \omega_j=
    \begin{cases}
        (\omega_j^q,0) & j=1,\dots,m \\
        (0,\omega_j^p) & j=m+1,\dots,2m
    \end{cases},
\end{equation*}
being $\omega_j^q,\omega_j^p\in\Vqp$ the Riesz representers of $\ell_j^q$ and $\ell_j^p$, respectively. We remark that this is not the most general setting, because of the assumption on the particular form of the measurements $\ell_j$. Nevertheless, we shall restrict our analysis to the case of separate measurements of $u^q$ and $u^p$ as this choice reflects the product space structure of $V$ and ultimately results in a simple expression of the gradient of $\beta^2$. With these definitions, we can write the Gram matrix $\bB=\bG(\Wtwom,\Vtwon)\in\bR^{2m\times2n}$ as
\begin{equation*}
    \bB = \begin{bmatrix}\bB_q \\ \bB_p\end{bmatrix} \qquad \text{ where } \qquad (\bB_q)_{i,j} := \inprodV{\omega_i^q}{v_j^q}_{\Vqp} \qquad (\bB_p)_{i,j}:=\inprodV{\omega_{i+m}^p}{v_j^p}_{\Vqp}
\end{equation*}
for $i=1,\dots,m$ and $j=1,\dots,2n$. Analogously, $\bA=\bG(\Wtwom,\Wtwom)$ is of the form
\begin{equation*}
    \bA = \begin{bmatrix}\bA_q & 0 \\0 & \bA_p\end{bmatrix} \qquad \text{ where } \qquad (\bA_q)_{i,j} := \inprodV{\omega^q_i}{\omega^q_j}_{\Vqp} \qquad (\bA_p)_{i,j}:=\inprodV{\omega^p_{i+m}}{\omega^p_{j+m}}_{\Vqp},
\end{equation*}
for $i,j=1,\dots,m$. We remark that the block diagonal form of $\bA$ is a consequence of the choice of the measurement $\ell_i$ as in \eqref{eq:liHam}. Therefore, we obtain
\begin{equation*}
    \bM(\Vtwon,\Wtwom) = \bB^\top\bA^{-1}\bB = \bB_q^\top\bA_q^{-1}\bB_q+\bB_p^\top\bA_p^{-1}\bB_p=:\bM_q+\bM_p\in\bR^{2n\times2n}.
\end{equation*}
The expression of $\nabla_\xs\beta^2$ can now be obtained using a similar derivation as in \Cref{sec:dyn_sens_plac}: if $c\in\bR^{2n}$ is the unit eigenvector corresponding to the smallest eigenvalue of $\bM(\Vtwon,\Wtwom)$, then for $\ell=1,\dots,d$ we define
\begin{equation*}
    \nabla_{\xs^{(\ell)}}\beta^2(\xs):= 2(\bA_q^{-1}\bB_qc)\odot(\bB^{(\ell)}_{q,D}c-\bA^{(\ell)}_{q,D}\bA_q^{-1}\bB_qc)+2(\bA_p^{-1}\bB_pc)\odot(\bB^{(\ell)}_{p,D}c-\bA^{(\ell)}_{p,D}\bA_p^{-1}\bB_pc)
\end{equation*}
where
\begin{equation*}
    (\bA^{(\ell)}_{q,D})_{i,j} := \inprodV{\frac{\partial\omega_i^q}{\partial \xs_i^{(\ell)}}}{\omega^q_j}_{\Vqp} \qquad (\bB^{(\ell)}_{q,D})_{i,j} := \inprodV{\frac{\partial\omega_i^q}{\partial \xs_i^{(\ell)}}}{v_j}_{\Vqp}.
\end{equation*}
The matrices with subscript $p$ are similarly defined by replacing $q$ with $p$. Finally, the gradient $\nabla_\xs\beta^2\in\bR^{2md}$ is defined as in \eqref{eq:gradbetasq}.

\bibliographystyle{plain}
\bibliography{references}

\end{document}